\documentclass[a4paper,oneside,12pt]{article}
\usepackage{amsmath,amsfonts,amssymb,amsthm,mathrsfs}
\usepackage[a4paper,vmargin={3.5cm,3.5cm},hmargin={2.5cm,2.5cm}]{geometry}
\usepackage[font=sf, labelfont={sf,bf}, margin=1cm]{caption}
\usepackage{graphicx,graphics}
\usepackage{epsfig}
\usepackage{latexsym}
\usepackage[applemac]{inputenc}
\usepackage{ae,aecompl}
\usepackage[english]{babel}
 \usepackage[colorlinks=true]{hyperref}
\usepackage{pstricks}
\usepackage{enumerate}

\newtheorem*{theorem}{Theorem}
\newtheorem{thm}{Theorem}

\newtheorem{proposition}[thm]{Proposition}
\newtheorem{lemma}[thm]{Lemma}
\newtheorem{cor}[thm]{Corollary}

\def\llbracket{[\hspace{-.10em} [ }
\def\rrbracket{ ] \hspace{-.10em}]}
\newtheorem{conjecture}{Conjecture}
\newtheorem*{rek}{Remark}

\newcommand{\td}{\tilde{\gamma}}
\newcommand{\op}[1]{\operatorname{#1 }}

\date{}

\title{\bf Simple random walk on the uniform infinite planar quadrangulation: Subdiffusivity \emph{via} pioneer points}
\begin{document}

\author{Itai Benjamini\footnote{Weizmann Institute of Science} \ and Nicolas Curien
\footnote{École Normale Supérieure}}
\maketitle

\abstract{We study the pioneer points of the simple random walk on the uniform infinite planar quadrangulation (UIPQ) using an adaptation of the peeling procedure of \cite{Ang03} to the quadrangulation case. Our main result is that, up to polylogarithmic factors, $n^3$ pioneer points have been discovered before the walk exits the ball of radius $n$ in the UIPQ. As a result we verify the KPZ relation \cite{KPZ88} in the particular case of the pioneer exponent and prove that the walk is subdiffusive with exponent less than $1/3$. Along the way, new geometric controls on the UIPQ are established.}
\section*{Introduction}

The goal of this work is to study the simple random walk on large random planar maps and especially on the Uniform Infinite Planar Quadrangulation (UIPQ). We show that the walk is dramatically affected by the geometry of the underlying random lattice and exhibits a  behavior very different from the classical deterministic Euclidean setting. For example, we show that the walk is \emph{subdiffusive}. Let us start by recalling the definition of the UIPQ.

A planar map is a proper embedding of 
a finite connected graph in the
two-dimensional sphere, considered up to orientation-preserving
homeomorphisms of the sphere. 
A \emph{quadrangulation} is a planar map whose faces all have degree $4$ (with the convention that 
if an edge lies entirely into a face then this edge is counted twice 
in the degree of the face). In  this work, the maps that we will consider will systematically be rooted, that is, given with a distinguished oriented edge $ \vec{e}$ called the root of the map. The origin vertex of the root edge is called the origin of the map and is denoted by $\rho$. 

The mathematical theory of random planar maps has been considerably growing over the last years motivated by the physics theory of 2D quantum gravity \cite{ADJ97}. In particular, Miermont and Le Gall recently proved that a large class of random planar maps properly rescaled converge towards a universal continuous random surface called the \emph{Brownian Map} \cite{LG11,Mie11}.  In this work, we choose a different perspective and study \emph{local limits} of random maps as introduced in \cite{BS01}.  If $m,m'$ are two rooted maps, the local distance between $m$ and $m'$ is 
 \begin{eqnarray*} \mathrm{d_{map}}(m,m') &=& \Big(1+ \sup\{r\geq 1 : \mathrm{Ball}(m,r) = \mathrm{Ball}(m',r)\}\Big)^{-1}, \end{eqnarray*}
where $ \mathrm{Ball}(m,r)$ denotes the map formed by the faces of $m$ that have at least one vertex at distance strictly less than  $r$ from the origin $\rho$ of $m$. The set of all finite quadrangulations is not complete for the metric $ \mathrm{d_{map}}$ and we have to add infinite quadrangulations to make it complete, see \cite{CMMinfini} for more details. Let $Q_n$ be a random rooted quadrangulation uniformly distributed over the finite set of all rooted quadrangulations with $n$ faces. Krikun \cite{Kri05} proved the following convergence  in distribution in the sense of $ \mathrm{d_{map}}$
 \begin{eqnarray} Q_n &\xrightarrow[n\to \infty]{(d)}& Q_\infty, \label{def:UIPQ}\end{eqnarray} where  $Q_\infty$  is a random infinite rooted quadrangulation called the \emph{Uniform Infinite Planar Quadrangulation} (UIPQ). See also the pioneer work of Angel \& Schramm \cite{AS03} who introduced a similar object (the UIPT) in the triangulation case. It is believed that the UIPT and the UIPQ share the same large-scale properties. However, we chose to focus on the UIPQ rather than on the UIPT because of the existence of ``nice" bijections between quadrangulations and simpler objects such as labeled trees \cite{Sch98} (these bijections do exist in the triangulation case but are less easy to manipulate). In particular, after the initial approach of Krikun \cite{Kri05}, Chassaing \& Durhuus \cite{CD06} gave a Schaeffer-like construction of the UIPQ based on a random infinite tree with positive labels (which was shown to be equivalent to that of Krikun in \cite{Men08}). The positivity constraint on the labels was relieved in \cite{CMMinfini} yielding to a third construction of the UIPQ (see Section \ref{schaeffer}). \medskip 


 The geometry of the UIPQ is very intriguing and is not completely understood. For instance, the UIPQ has a striking  growth rate of $r^4$ \cite{CD06,LGM10} but yet possesses separating cycles of linear length at all scales \cite{Kri05}.  These isoperimetric inequalities  heuristically suggest that the UIPQ has many folds and bottlenecks at all scales in which the nearest-neighbor simple random walk (SRW) could be trapped for a while, slowing it down. We will study this slowing effect by looking at particular points of the SRW called \emph{pioneer points}. Let us define properly this notion.
 
  Conditionally on $Q_\infty$, let $(X_{n})_{n\geq 0}$ be a nearest-neighbor simple random walk on $Q_{\infty}$ starting from the origin $\rho$. For any $k \geq 0$ we denote by $ \mathcal{R}_k$ the set of all faces of $Q_\infty$ that are adjacent to  the range $\{X_0,X_1, \ldots , X_k\}$ of the walk up to time $k$. A time $k\geq 1$ is a pioneer time (in which case we say that $X_k$ is a pioneer point) if $X_k$ lies on the boundary of the only infinite component of $Q_\infty \backslash \mathcal{R}_{k-1}$ (the UIPQ has almost surely one end \cite{Kri05}). Our main result is:
 
\begin{thm}[Main result]\label{main} Let $Q_{\infty}$ be the uniform infinite planar quadrangulation.  Conditionally on $Q_{\infty}$, let $(X_{n})_{n\geq 0}$ be a nearest-neighbor simple random walk on $Q_{\infty}$ starting from $\rho$. We denote by $P_1, P_2, \ldots$ the pioneer points of $(X_n)_{n \geq 0}$. Then there exists a constant $\kappa >0$ such that a.s. we eventually  have  \begin{eqnarray*} n^{1/3} \log^{-\kappa}(n) \quad \leq \quad \max_{0 \leq k \leq n} \mathrm{d_{gr}}(\rho, P_k)  \quad \leq \quad n^{1/3}\log^\kappa(n).  \end{eqnarray*} \end{thm}

We did not try to compute the best value of $\kappa$ given by our proof and we do not have a precise guess for the correct logarithmic fluctuations. As a corollary of the proof of Theorem \ref{main} we have:

\begin{cor}[Subdiffusivity] \label{subdiffusive}With the notation of Theorem \ref{main}, there exists a constant $\kappa'>0$ such that a.s. we eventually have
   \begin{eqnarray*} \label{subdiff} \op{d^{
 }_{gr}}(\rho,X_{n})  & \leq &\quad n^{1/3} \log^{\kappa'}(n).  \end{eqnarray*}
\end{cor}



The simple random walk on the UIPQ thus has a \emph{subdiffusive} behavior since it displaces much slower than the ${n}^{1/2}$ classical behavior of the simple random walk on $ \mathbb{Z}^d, d \geq 1$. This phenomenon has first been suggested by Pierre-Gilles De Gennes \cite{DG76} for the simple random walk on a critical percolation cluster:  \emph{``la fourmi dans un labyrinthe"}. This was rigorously proved by Kesten \cite{Kes86} for simple random walk on the infinite incipient cluster of critical two-dimensional Euclidean Bernoulli percolation (the exact exponent is still unknown) and finite variance critical Galton-Watson trees conditioned to survive (exponent ${1/3}$). This phenomenon has then been established for others models, see e.g.\,\cite{Bar04,BK06,CK08,K10}. We do not expect the $1/3$ exponent of Corollary \ref{subdiffusive} to be sharp and conjecture that $1/4$ is the correct value: \begin{conjecture} \label{1/4}The subdiffusivity exponent of the SRW on the UIPQ is $1/4$:   \begin{eqnarray*} \max_{0 \leq k \leq n}\mathrm{d_{gr}}(\rho,X_{k}) &\approx &n^{1/4}.  \end{eqnarray*}
\end{conjecture}See Section \ref{comments} for comments. \\

 Usually, the road map to prove a subdiffusivity result is to estimate the volume growth and  resistances in the graph. In our setting, evaluating  resistances in the UIPQ remains a challenging  problem. In particular, it is still open to show that the resistance between $\rho$ and $\infty$ is infinite, in other words:
 \begin{conjecture} \label{recurrence}\cite{AS03}  The UIPQ is almost surely recurrent.
\end{conjecture}
 Rather than estimating resistances, the key to prove Theorem \ref{main} it to use one of the main features of random planar maps: the \emph{spatial Markov property} (Theorem \ref{markov}). This property can roughly be stated as follows: Imagine that we explore a simply connected region of the UIPQ, then conditionally on the length of the boundary of this region, the remaining part of the UIPQ is independent of the explored region. The spatial Markov property of random planar triangulations has been used by Angel \cite{Ang03} to study several properties of the UIPT via the so-called \emph{peeling process}. This is a clever random algorithm that discovers step-by-step the UIPT by revealing one face at a time, like ``peeling an apple'' \cite{Wat95}. Using the explicit transition probabilities of the peeling process of the UIPT, Angel \cite{Ang03} has obtained sharp estimates (up to polylogarithmic fluctuations) of the perimeter and the size of the triangulation discovered after $n$ steps of peeling.

  In this work, we will develop the same approach in the quadrangulation case and provide the analogs of the peeling estimates of Angel in the case of the UIPQ (Theorem \ref{estimpeeling}). However, our tactics here does not consist in mimicking  the proofs of \cite{Ang03} but rather to use the universality of the peeling process in order to translate geometric controls on the UIPQ (Section \ref{geo}) into estimates on the peeling process.
 \bigskip

  Let us give a rough sketch of the proof of our main result. The idea is to discover the UIPQ along a simple random walk path using the peeling device. During this exploration, roughly speaking, only the pioneer points of the walk trigger the discovery of a new quadrangle. 
 It then turns out that the boundary of the quadrangulation discovered after peeling $n$ quadrangles of $Q_{\infty}$ (or equivalently, after discovering $\approx n$ pioneer points of the walk) is of order $n^{2/3}$ (see Theorem \ref{estimpeeling}). Now, by the spatial Markov property of the UIPQ, conditionally on the boundary of length $ \approx n^{2/3}$, the remaining part of $Q_{\infty}$ is independent of the revealed part. It has recently been  proved in \cite{CMboundary} that the typical distance between boundary points in a UIPQ with a boundary of perimeter $p=n^{2/3}$ is of order $\sqrt{p}=n^{1/3}$ which is the first glimpse at the $n^{1/3}$ appearing in Theorem \ref{main}.\bigskip

Random planar maps are key  tools in  understanding Euclidean statistical physics  systems $via$ the Quantum Gravity approach, see e.g. \cite{ADJ97}. Especially, the KPZ formula (Knizhnik,  Polyakov and Zamolodchikov \cite{KPZ88}) predicts relations between critical exponents of statistical mechanics models on a Euclidean lattice and the analogs on a random lattice (a random map).  Duplantier and Sheffiled rigorously proved  the KPZ relations in a random geometry constructed from the Gaussian free field \cite{DS09}. One missing link is the connection between random planar maps and the Gaussian free field, see the conjectures in \cite{DS09,She10} and \cite{Ben10}. We propose a verification of the KPZ formula concerning pioneer points exponent of the simple random walk (Section \ref{comments}). 
We also use the KPZ prediction on the disconnection exponent of the simple random walk on a random lattice to support Conjecture \ref{1/4}.  
\\ 

The paper is organized as follows. The next section introduces the spatial Markov property of the UIPQ, the peeling process, and the main estimates about it (Theorem \ref{estimpeeling}). Section 	\ref{schaeffer} presents the construction of \cite{CMMinfini} that we use in Section 	\ref{geo} to give new geometric lemm\ae\ on the UIPQ. In particular, we study the vertex degrees in the UIPQ  (Proposition \ref{degree}) and provide uniform control on the volume growth (Proposition \ref{growth}) and on the length of the separating cycles at a given height (Proposition \ref{krikun}). We then proceed to the (very short) proofs of Theorems \ref{estimpeeling} and \ref{main} in Section \ref{proofs}. Unsurprisingly, the final section contains comments (in particular about the KPZ relation) and open questions.\\

\noindent \textbf{Acknowledgments:} We are grateful to Christophe Garban for a very stimulating discussion and for useful remarks. Thanks also go to Jean-Fran\c cois Le Gall for precious comments on an early version of this work and to Igor Kortchemski for a careful reading.

\section{The peeling process}

\subsection{Spatial Markov property}
This section is adapted from \cite[Sections 4 and 5]{AS03}. Since the proof are \emph{mutatis mutandis} the same as in the triangulation case we leave them to the reader. Let us first introduce a few notions.

A quadrangulation with $k\geq 1$ holes is a (rooted) map with $k$ distinct distinguished faces called the holes of the map, and such that all the non-distinguished faces have degree four. In the following, we will always assume that the boundaries of the holes are cycles with no self-intersection. Notice that by bipartiteness the holes must be of even degree.   A quadrangulation $q$ with one hole is also called a quadrangulation with a (simple) boundary. 
In this case, the degree $|\partial q|$ of the unique hole is called the perimeter  of the map and its size $|q|$ is its number of inner faces. A quadrangulation with a simple boundary of perimeter $2p$ will also be called a quadrangulation of the $2p$-gon. By convention, all the quadrangulations with a boundary that we consider here are rooted on the boundary in such a way that the hole is lying on the right-hand side of the distinguished oriented edge $\vec{e}$. Note that a quadrangulation of the $2$-gon can be considered as a rooted quadrangulation (without hole) by contracting the unique face of degree $2$. 

\begin{figure}[!h]
 \begin{center}
 \includegraphics[width=10cm]{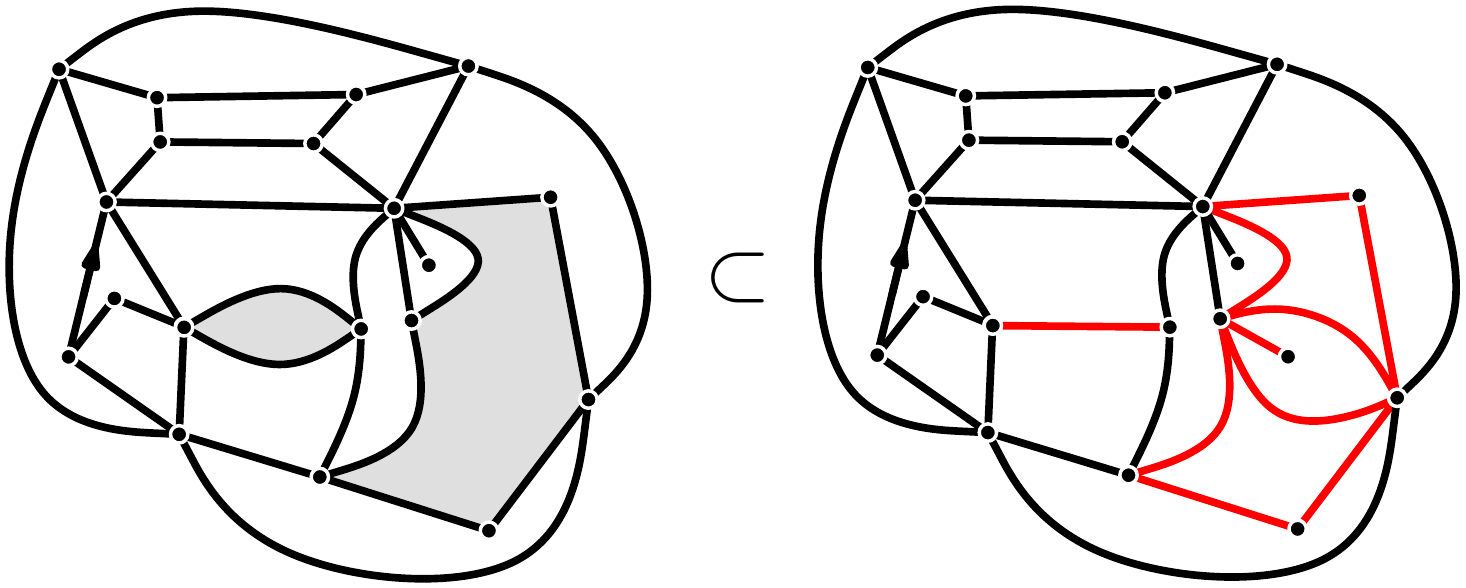}
 \caption{ \label{submap}A (rigid) quadrangulation with holes and a way of filling it. Remark that the $2$-gon on the left has been closed by a simple boundary of perimeter $2$ and size $0$.}
 \end{center}
 \end{figure}

\paragraph{Enumeration results.}  We write 	${\mathcal{Q}}_{n,2p}$ for the set of  all rooted quadrangulations of the $2p$-gon with $n$ inner faces. From \cite[(2.11)]{BG09} we read that for $p \geq 1$ and $n \geq 1$ we have  \begin{eqnarray} \label{enumeration} \# { \mathcal{Q}}_{n,2p} &=&  \frac{3^{n-p}(3p)!(2n+p-1)!}{p!(2p-1)!(n-p+1)!(n+2p)!}.  \end{eqnarray}
In the case $p=1$ and $n=0$, the only element of $\mathcal{Q}_{0,2}$ is the quadrangulation with a simple boundary composed of one oriented edge. It is not ``really'' a quadrangulation but will be interpreted as follows (see the remark after Proposition 1.6 in \cite{Ang03}): A quadrangulation  of the $2p$-gon will often be used to close a hole of degree $2p$ in another quadrangulation, thus in the case $p=1$ and $n=0$ we close a hole of length $2$ by gluing its two edges together. See Fig.\,\ref{submap}.

More precisely, if $q$ is a quadrangulation with holes we choose, once for all, a deterministic way of distinguishing an oriented edge on the boundary of each hole of $q$ (such that the hole lies on the left of it). If $q$ has a hole of perimeter $2p$ and if we are given a quadrangulation $h$ of the $2p$-gon,  we can glue $h$ inside the hole of $q$ by identifying their boundaries (such that the oriented edge of the hole coincides with that of $h$). See Fig.\,\ref{submap}. Let $q,Q$ two quadrangulations with holes. We say that  $q$  is a submap of $Q$ if $Q$ can be obtained by filling  some of the holes of $q$ with  quadrangulations with  simple boundaries and we write $$ q \subset Q.$$  We say that $q$ is \emph{rigid} if two different ways of filling it lead to two different planar maps (see \cite[Definition 4.7]{AS03}). An easy adaptation of \cite[Lemma 4.8]{AS03} yields that any (rooted) quadrangulation with a simple boundary is rigid.
 
From \eqref{enumeration}, the asymptotic of $\#  { \mathcal{Q}}_{n,2p}$ takes the form $\#  { \mathcal{Q}}_{n,2p} \sim  C_{2p} 12^n n^{-5/2}$ when $n \to \infty$, where $C_{2p}$ a positive constant depending only on $p$. The $n^{-5/2}$ polynomial correction is typical of the enumeration of planar maps and plays a crucial role in the large scale structure of the UIPQ. In particular, the series $$\sum_{n=0}^\infty \#  { \mathcal{Q}}_{n,2p}12^{-n},$$  is convergent and we denote its sum by $Z_{2p}$. 
Following \cite[Definition 2.3]{AS03} we define the \emph{free distribution} on rooted quadrangulations of the $2p$-gon as the probability measure $\nu_{2p}$ that assigns the weight  $ 12^{-n}{Z_{2p}^{-1}}$ to each element of $\cup_{n \geq 0} {\mathcal{Q}}_{n,2p}.$

The convergence \eqref{def:UIPQ} can easily be extended to the case of quadrangulations with simple boundary: If $Q_{n,2p}$ is a uniform quadrangulation with size $n$ and perimeter $2p$ then we have the convergence in distribution for $ \mathrm{d_{map}}$
 \begin{eqnarray*} Q_{n,2p} & \xrightarrow[n\to\infty]{(d)}& Q_{\infty,2p},  \end{eqnarray*} 
 where $Q_{\infty,2p}$ is the UIPQ with simple boundary of perimeter $2p$ of UIPQ or the $2p$-gon. This convergence is a simple consequence of \eqref{def:UIPQ}, see \cite{CMboundary}.  We can now state the spatial Markov property of the UIPQ. Recall that almost surely the UIPQ (and more generally the UIPQ of the $2p$-gon) has one end \cite{Kri05}.


\begin{thm}[Spatial Markov Property]\label{markov} Let $q$ be a finite rigid (rooted)  quadrangulation with $k+1$ holes of even degrees $p_{0},p_{1}, \ldots, p_{k}$ such that the $0$th hole is on the right of $ \vec{e}$. Conditionally on the event  $\{q \subset Q_{\infty,p_{0}}\}$ denote $H_{1}, H_{2}, \ldots , H_{k}$ the quadrangulations filling the first, second, third\ldots holes of $q$ in $Q_{\infty,p_{0}}$. Then conditionally on $\{ H_{i} \mbox{ is infinite}\}$,  the  quadrangulations $H_{1}, \ldots ,H_{k}$ are independent  and
\begin{enumerate}[(i)]
\item $H_{i}$ has the same distribution as $Q_{\infty,p_{i}},$
\item for $j\ne i$, $H_{j}$ is distributed according to $\nu_{p_{j}}$. \end{enumerate}
\end{thm}
\subsection{The peeling algorithm}
A growth algorithm for random maps, \emph{the peeling process}, was first used heuristically by physicists (see \cite{Wat95} and \cite[Section 4.7]{ADJ97}). Angel \cite{Ang05,Ang03} then  defined it rigorously and used it to study the volume growth and  site percolation on the uniform infinite planar triangulation. We adapt his ideas to the context of the UIPQ.\medskip 

The peeling process  is a procedure that allows us to discover step-by-step the UIPQ by revealing one quadrangle at a time. Formally, we construct (on the same probability space) the uniform infinite planar quadrangulation $Q_{\infty}$ and a sequence of  rooted quadrangulations with  a boundary $Q_{0}\subset Q_{1} \subset \ldots \subset Q_{n} \subset \ldots \subset Q_{\infty}$, such that for every $i \geq 0$, conditionally on $Q_{0}, \ldots , Q_{i}$, the remaining part $Q_{\infty}\backslash Q_{i}$  has the same distribution as $Q_{\infty,|\partial Q_{i}|}$. The sequence is  constructed inductively: \medskip 

The quadrangulation $Q_{0}$ is the root edge of $Q_{\infty}$, which can be viewed as a quadrangulation with a boundary of perimeter $2$. We write $\mathcal{F}_{n}$ for the filtration generated by $Q_{0}, \ldots , Q_{n}$. By the induction hypothesis, conditionally on $\mathcal{F}_{n}$, the remaining part $Q_{\infty}\backslash Q_{n}$  which is contained in the unique hole of $Q_{n}$ has the same distribution as $Q_{\infty,|\partial Q_{n}|}$. 

The conditional distribution of $Q_{n+1}$ knowing $\mathcal{F}_{n}$ can be described as follows. We first   choose deterministically, or with the help of a randomized algorithm  independent of $Q_{\infty} \backslash Q_{n}$, an edge $e^*$ on $ \partial Q_{n}$ and re-root $Q_{\infty} \backslash Q_{n}$ at this edge. Since the choice of $e^*$ is independent of $Q_{\infty}\backslash Q_{n}$ the newly rooted map still has the law of a UIPQ of the $|\partial Q_{n}|$-gon. The peeling process then reveals {the quadrangle in the remaining part $Q_{\infty}\backslash Q_{n}$ containing the edge $e^*$}. Three cases may happen (see Fig.\,\ref{cases}):

\begin{figure}[h]
\begin{center}
\includegraphics[width=10cm]{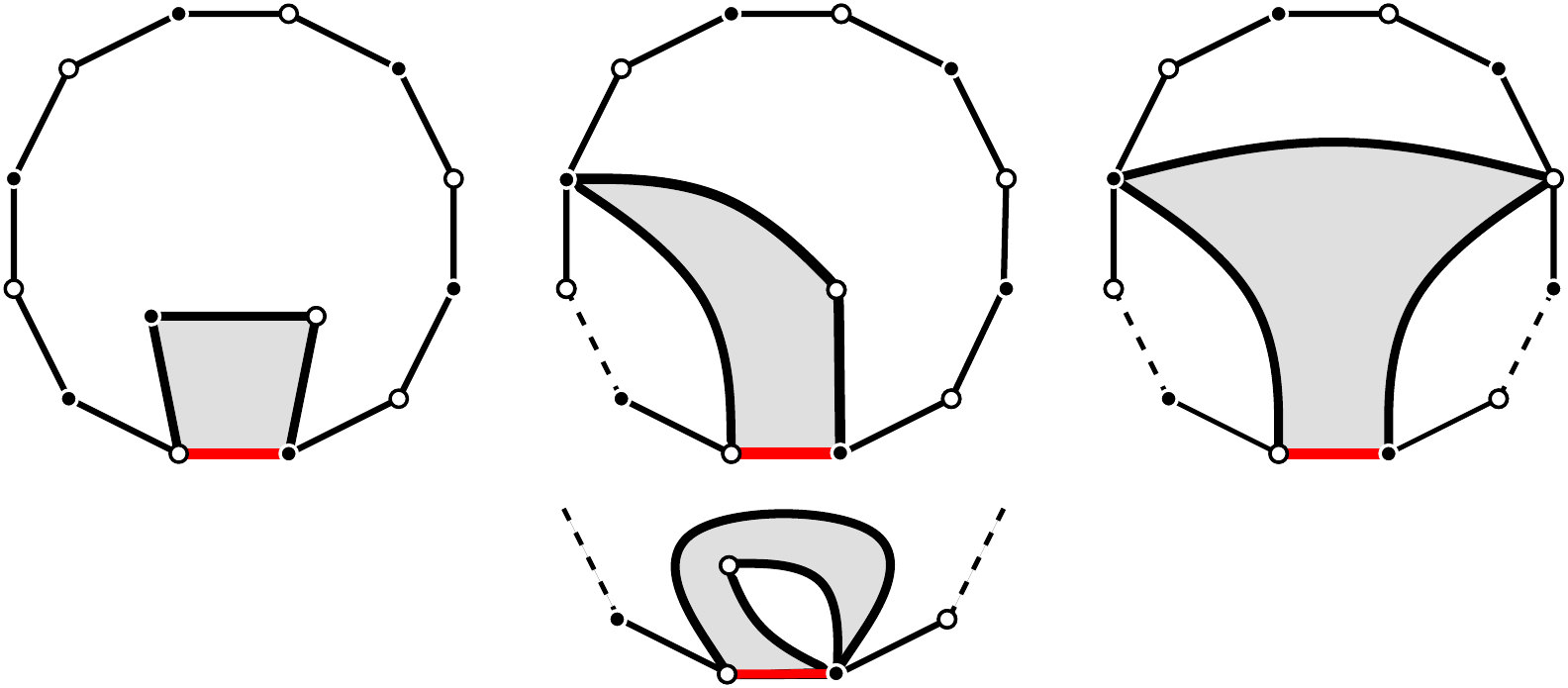}
\end{center}
\caption{ \label{cases}Illustration of the three cases that could happen when we reveal a new quadrangle in the unknown region. Notice that in the second case, two of the points on the boundary may be confounded.}
\end{figure}

\begin{itemize}
\item The revealed quadrangle has two vertices lying on the boundary $\partial Q_{n}$.  In this case, we set $Q_{n+1}$ to be the union of $Q_{n}$ together with the revealed quadrangle. Hence  $Q_{n+1}$ is a quadrangulation with a boundary of length $|\partial Q_{n}|+2$, and thanks to Theorem \ref{markov}, conditionally on this event and on $\mathcal{F}_{n}$, the remaining quadrangulation $Q_{\infty} \backslash Q_{n+1}$ has the same distribution as $Q_{\infty,|\partial Q_{n+1}|}$.
\item The quadrangle has three vertices lying on the boundary (two of these vertices might be identified). This quadrangle thus separates the remaining part into two quadrangulations $Q_{n,1}$ and $Q_{n,2}$ which are respectively quadrangulations of the $p_{1}$-gon and $p_{2}$-gon, such that $p_{1}+p_{2} = |\partial Q_{n}|+2$. Since $Q_{\infty}$ almost surely has one end, only one of this two components is infinite. For definiteness we argue on the event $$A=\{ Q_{n,1} \mbox{ is finite}, Q_{n,2} \mbox{ is infinite}\}.$$ Thanks to Theorem \ref{markov}, conditionally on the revealed quadrangle, on $A$ and on $\mathcal{F}_{n}$, $Q_{n,1}$  is distributed according $\nu_{p_{1}}$ and is independent of $Q_{n,2}$ which has the same distribution as $Q_{\infty,p_{2}}$. We thus set  $Q_{n+1}$ to be the union of $Q_{n}$, of  $Q_{n,1}$ and of the revealed quadrangle. Notice that $Q_{n+1}$ is a quadrangulation with a boundary of perimeter $p_{2}$ and that conditionally on  $\mathcal{F}_{n+1}$, $Q_{\infty} \backslash Q_{n+1}$ has the same distribution as $Q_{\infty,|\partial Q_{n+1}|}$.
\item The quadrangle has its four vertices lying on the boundary of $Q_{n}$ and separates the remaining part into three quadrangulations $Q_{n,1}$, $Q_{n,2}$ and $Q_{n,3}$. Similarly as in the preceding case, only one of these quadrangulations is infinite.  We then set $Q_{n+1}$ to be the union of these finite quadrangulations and of the revealed quadrangle and check that $Q_{\infty} \backslash Q_{n+1}$ has the desired law.
\end{itemize}

\indent We stress the fact that there are many ways to do the peeling  of $Q_{\infty}$ according to the algorithm we use to choose the next quadrangle to reveal (provided that this choice is independent of the unknown part $Q_{\infty} \backslash Q_{n}$). Although the distribution of $Q_{0}, \ldots , Q_{n}\ldots $ may depend on the algorithm,  the process $(|\partial Q_{n}|)_{n\geq 0}$ is actually a Markov chain whose distribution does not depend on the manner we revealed the squares in $Q_{\infty}$. Moreover the volume $|Q_{n}|$ of $Q_{n}$ (its number of vertices) is obtained from $|Q_{n-1}|$ by filling with free quadrangulations of proper perimeters (and independent of the past) one or two holes whose perimeter only depend on the quadrangle revealed at time $n-1$. Therefore a moment's thought shows that  $(|\partial Q_{n}|,|Q_{n}|)_{n\geq 0}$ is a homogeneous Markov chain whose transition probabilities do not depend on the procedure chosen to do the peeling. Thus we have:

\begin{lemma}\label{peelingloi} For any peeling procedure $Q_{0}, \ldots , Q_{n}, \ldots$ the process $(|\partial Q_{n}|, |Q_{n}|)_{n\geq 0}$ has the same distribution. \end{lemma}

In \cite{Ang03}, Angel explicitly computed the transition probabilities of the peeling in the case of the UIPT. Through a careful analysis of this chain, he proved that the boundary and the size of the triangulation obtained after $n$ steps of peeling are respectively of order $ n^{2/3}$ and $n^{4/3}$ up to polylogarithmic fluctuations. We will prove the same result  in the case of the UIPQ. Before that, let us acquaint the reader with a useful notation. \medskip

\label{notations}
In all this paper if $(Y_{n})_{n \geq 0}$ is a  random process indexed by $ \mathbb{N}$ with values in $ \mathbb{R}_{+}$, we write $Y_{n} \succeq n^\alpha$ resp.  $Y_{n} \preceq n^\alpha$ for $\alpha >0$ if there exists a constant $\kappa >0$ such that we almost surely have    \begin{eqnarray*}\lim_{n \to \infty} \frac{Y_{n}}{n^{\alpha}\log^{-\kappa}(n)} = \infty \quad  \mbox{resp.} \quad  \lim_{n \to \infty} \frac{Y_{n}}{n^{\alpha}\log^{\kappa}(n)} = 0.  \end{eqnarray*} If we have both $ Y_{n} \preceq n^\alpha$ and $n ^\alpha \preceq Y_{n}$ we write $Y_{n} \approx n^\alpha$. In words, $Y_n \approx n^{\alpha}$ means that almost surely $Y_n$ grows like $n^\alpha$ up to polylogarithmic fluctuations.  We also recall the classical Landau notation $x_n = O(y_n)$ (resp. $x_n = \Theta(y_n)$) if there exists a constant $0<C<\infty$ (resp.\,$0<c<C<\infty$) such that $x_n \leq Cy_n$ (resp.\,$cy_n \leq x_n \leq C y_n$). We also denote $x_{n} \sim y_{n}$ if the quotient $x_{n}/y_{n}$ goes to $1$ as $n \to \infty$.
\begin{thm} \label{estimpeeling}For any peeling procedure $Q_{0}, \ldots , Q_{n}, \ldots$ we have  \begin{eqnarray} \label{estim}|\partial Q_{n}| \approx n^{2/3} \qquad \mbox{and} \qquad  |Q_{n}| \approx n^{4/3}. \end{eqnarray}\end{thm}

Using the enumeration results of \cite{BG09} it is possible to explicitly compute  the transitions probabilities of the Markov chain $ |\partial Q_{n}|$ (see \cite{ACpercopeel}). It is also believable that the arguments of \cite{Ang03} could be adapted to show Theorem \ref{estimpeeling}. However, this is not the path we are about to follow.  We propose a softer approach to this result. The idea is to get estimates on the peeling process \emph{via} geometric estimates on the UIPQ. Indeed, because of Lemma \ref{peelingloi} it is sufficient to prove  Theorem  \ref{estimpeeling} for \emph{one} carefully chosen peeling algorithm.  We will thus introduce an adaptation of the method proposed in \cite{Ang03} to analyze the volume growth in the UIPT: After establishing new results  on the volume growth in $Q_{\infty}$, this process will be used in Section \ref{proofs} 	 to prove Theorem \ref{estimpeeling}. These estimates will then be used with a second peeling  coupled with a simple random walk  on $Q_{\infty}$.

\subsection{Peeling by layers} \label{peelL}In this section we present a peeling process that discovers $Q_{\infty}$ ``layer after layer''. It is an adaptation of the peeling procedure of \cite[Section 2]{Ang03}.  Together with the geometric estimates of Section \ref{geo}, this process will be used to deduce Theorem \ref{estimpeeling}. In order to describe this peeling, we just have to tell how do we choose the next quadrangle to reveal in the process. We will then interpret it in a more geometric way.
 
 \paragraph{Algorithm $\mathcal{L}$.}
\begin{quote}\underline{Algorithm $\mathcal{L}$ ``Layer''}: Assume that $Q_{n} \subset Q_{\infty}$ is the quadrangulation with a boundary $\partial Q_{n}$ containing the root edge of $Q_{\infty}$ given by the peeling procedure $ \mathcal{L}$ at time $n$. The next quadrangle to reveal is chosen as follows. Pick an edge $e^*$ on the boundary $\partial Q_{n}$ such that one of its extremity minimizes $\{ \op{d_{gr}}(\rho, x) : x \in \partial Q_{n}\}$ and reveal the quadrangle in $Q_{\infty}\backslash Q_{n}$ that contains $e^*$. Notice that there might be several edges satisfying this property, in this case, choose deterministically one of them. \end{quote}

This algorithm thus gives a way to peel the UIPQ. However, one must be careful with this procedure since one could have to use the ``unknown'' part $Q_{\infty}\backslash Q_{n}$ in order to compute the graph distance between the origin $\rho$ and a point $x$ on the boundary  $\partial Q_{n}$. Recall that the choice of the edge to peel must not depend on $Q_{\infty} \backslash Q_{n}$  otherwise the law of the sequence $ ( |\partial Q_n|, |Q_n|)_{n \geq 0}$ might not be the same as a standard peeling process (Lemma \ref{peelingloi}). Fortunately, we will see (Proposition \ref{ouff} $(ii)$) that if the preceding algorithm \emph{has been used from the very first step $n=0$}, then the graph distance between any point on $\partial Q_n$ to $\rho$ can be computed using $Q_n$ only and thus the preceding algorithm yields to a true peeling process. Let us first introduce a piece of notation. 

\paragraph{Interpretation.} Recall that for $r \geq 1$, we denote by $\op{Ball}(Q_\infty,r)$ the quadrangulation contained in $Q_{\infty}$ composed of the faces that have at least one vertex at distance strictly less than $r$ from the origin $ \rho \in Q_{\infty}$. In particular $ \mathrm{Ball}(Q_\infty,r) \subset \{ u \in Q_\infty : \mathrm{d_{gr}}(\rho,u) \leq r+1\}$ in terms of vertex sets. By convention, the root edge of $Q_{\infty}$ is considered as a face of degree two, so that $ \vec{e} \subset \mathrm{Ball}(Q_{\infty},1)$. 

A moment's thought shows that $\op{Ball}(Q_{\infty},r)$ is a quadrangulation with  holes, in particular the boundaries of the holes are cycles with no self-intersection. Since $Q_{\infty}$ almost surely has one end,  only one hole of $\op{Ball}(Q_{\infty},r)$ corresponds to an infinite quadrangulation of $Q_{\infty}\backslash \op{Ball}(Q_{\infty},r)$. We denote the boundary of this hole by $\gamma_{r}$ and  called it the \emph{separating cycle of $ \rho$ and $\infty$ in $Q_{\infty}$ at height $r$}. Observe that this cycle is actually a simple path that alternatively visits  vertices at distance $r$ and $r+1$ from the origin. We also denote by $\overline{\op{Ball}}(Q_\infty,r)$ the quadrangulation obtained from $\op{Ball}(Q_{\infty},r)$ by filling all the finite holes of $ \mathrm{Ball}(Q_{ \infty},r)$ with their respective quadrangulations in $Q_{\infty}$. We call $ \overline{ \mathrm{Ball}}(Q_{\infty},r)$ the \emph{hull} of the ball of radius $r$ in $Q_{\infty}$. See Fig.\,\ref{boules} below.

\begin{figure}[h]
\begin{center}\includegraphics[width=15cm]{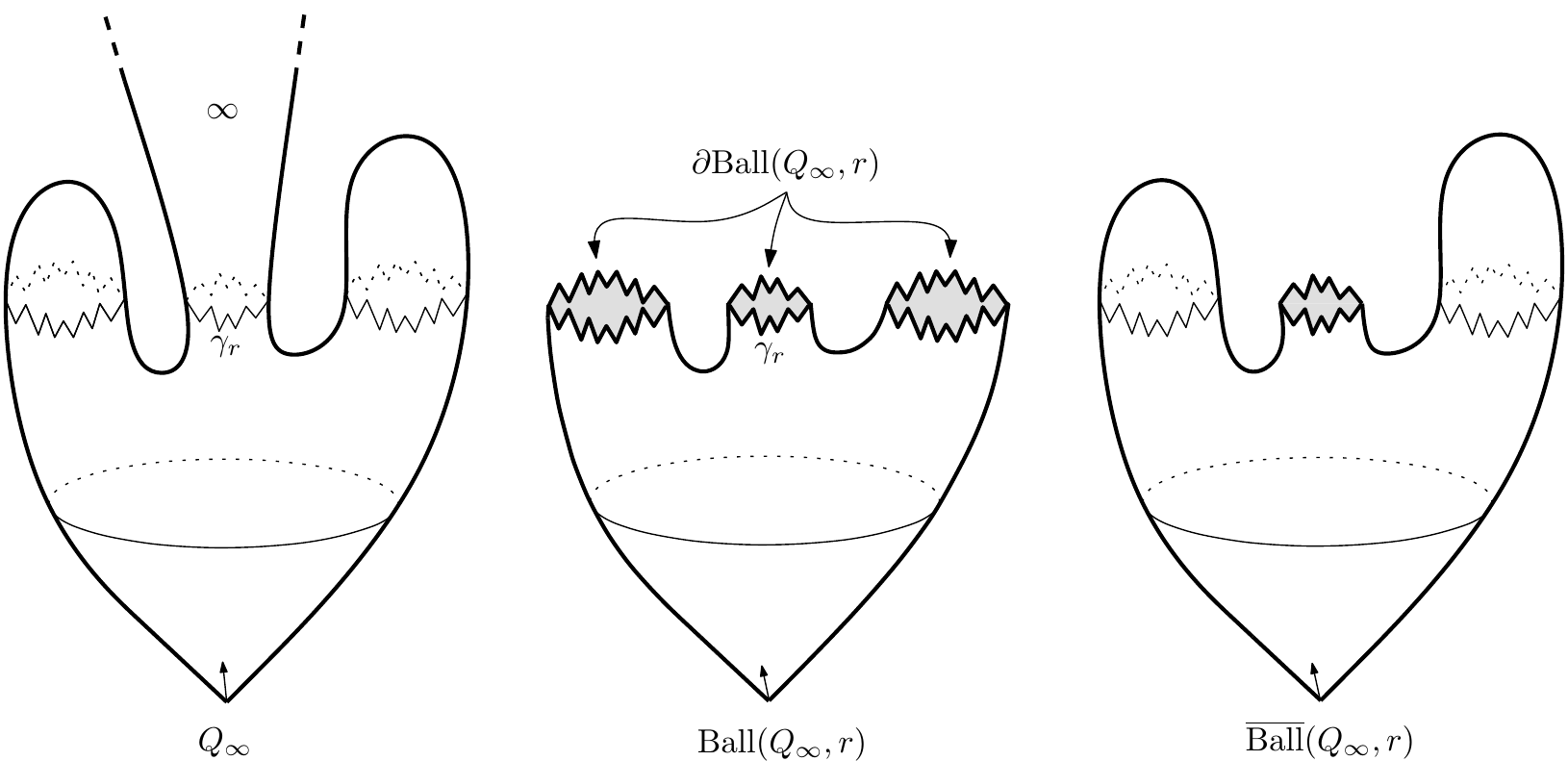}
\caption{ \label{boules}Illustration of $\op{Ball}(Q_{\infty},r)$, $\gamma_{r}$, and $\overline{\op{Ball}}(Q_{\infty},r)$ in $Q_\infty$.}
\end{center}
 \end{figure}
 The peeling process under Algorithm $ \mathcal{L}$ can geometrically be  interpreted as follows: It roughly discovers $Q_{\infty}$ layer after layer and stays very close to the cycles $\gamma_{k}$, $k \geq 0$ (see the figures of Section 4.7 in \cite{ADJ97}). More precisely:

\begin{proposition} \label{ouff} Let $Q_0 \subset Q_1 \subset Q_2 \subset \ldots \subset Q_\infty$ be the successive quadrangulations with a boundary discovered using Algorithm $ \mathcal{L}$. Then we have:

\begin{enumerate}[(i)]
\item  For every integer $r\geq 1$, let  $T_{r}$ be the first $n\geq 0$ such that  $\op{d_{gr}}( \rho,\partial Q_{n}) \geq  r$. Then $T_r <\infty$ a.s. and we have
\begin{eqnarray}
\label{cycleinf} Q_{T_r} \quad = \quad \overline{\mathrm{Ball}}(Q_\infty,r) \qquad \mbox{ and } \qquad \partial Q_{T_{r}} \quad = \quad  \gamma_{r}. \nonumber 
\end{eqnarray}
\item Futhermore, for any $n \geq 0$ and for any $x \in \partial Q_n$, the graph distance $ \mathrm{d_{gr}}(\rho,x)$ is measurable with respect to $Q_n$ and thus the sequence $(|\partial Q_n|, |Q_n|)_{n \geq 0}$ has the distribution of a standard peeling process.
\end{enumerate}
 \end{proposition}
\proof[Proof (Sketch)] We prove the proposition by induction on $r \geq 1$. Let us first examine the case $r =1$. We have $Q_0 = \vec{e}$ and start discovering some of the quadrangles that contain $\rho$. Notice that for any quadrangle adjacent to $ \rho$, the graph distance from $\rho$ of its vertices is either $\{ 0,1,0,1\}$ or $\{0,1,2,1\}$ by bipartiteness. Thus, when we discover a quadrangle containing $\rho$ one can deduce the graph distance from $\rho$ of its vertices by just looking at the quadrangulation discovered so far. Hence  as long as $n < T_1$ the graph distances of vertices of $\partial Q_n$ to $\rho$ are measurable with respect to $Q_n$. 

Furthermore, all the quadrangles discovered for $n < T_1$ as well as the holes they created (which are filled-in during the process) are contained in $ \overline{\mathrm{Ball}}(Q_\infty,1)$. We  stop at $T_{1}$ when the origin $\rho$ is not on the boundary of the current discovered quadrangulation $Q_{T_1}$.  By the remarks above we have $ Q_{T_1} \subset \overline{ \mathrm{ Ball}}(Q_\infty,1)$. The converse inclusion is deduced from the fact that the boundary of $Q_{T_1}$ is composed of vertices that are alternatively at distance $1$ and $2$ from $\rho$. Hence $Q_{T_1} = \overline{ \mathrm{ Ball}}(Q_\infty,1)$ and $\gamma_1 = \partial Q_{T_1}$.\\
The general case $r \geq 2$ is pretty much the same and is safely left to the reader. \endproof 

 It follows from Proposition \ref{ouff} that if $r$ is the minimal distance in $Q_\infty$ from a vertex in $\partial Q_{n}$ to the origin $\rho$ then we have \begin{eqnarray}  \op{d_{gr}}( \rho,x) &\in& \{r,r+1,r+2\}, \label{rr+1r+2}  \qquad \mbox{and} \qquad \op{d}_{\op{gr}}(\partial Q_{n}, \gamma_{r}) \leq
 2.\end{eqnarray}
\subsection{Peeling along a simple random walk} \label{peelW}We now describe another way of peeling  $Q_{\infty}$. This one is coupled with a simple random walk and discovers the quadrangulation when  necessary for the walk to make one more step. This peeling process is one of the keys in the proof of Theorem \ref{main}. We start with the formal definition of this algorithm and then interpret it in terms of pioneer points.

\paragraph{Algorithm $\mathcal{W}$. }
\begin{quote}\underline{Algorithm  $\mathcal{W}$ ``Walk'':} Let $Q_{\infty}$ be the uniform infinite planar quadrangulation and conditionally on $Q_{\infty}$, let $(X_{n})_{n\geq 0}$ be a nearest-neighbor simple random walk on $Q_{\infty}$ starting from $ \rho$. We do the peeling process each time we need it for the SRW to displace. More precisely, we define a sequence $$ \vec{e} = Q_{0} \subset Q_{1} \subset \ldots \subset Q_{n} \subset\ldots \subset Q_{\infty}$$ of quadrangulations with boundaries and two random non decreasing functions $f,g : \mathbb{N} \to \mathbb{N}$ such that $f(0)=g(0)=0$, $X_{g(k)} \in Q_{f(k)}$ for every $k\geq 0$, and whose evolution is described by induction as follows. 

We have two cases. If the current position $X_{g(k)}$ of the simple random walk belongs to $\partial Q_{f(k)}$, then choose an edge $e^*$ on $\partial Q_{f(k)}$ containing $X_{g(k)}$ and set $f(k+1):=f(k)+1$ and $g(k+1):=g(k).$ The quadrangulation $Q_{f(k+1)}$ is the map obtained after the peeling associated with the edge $e^*$. If the current position $X_{g(k)}$ of the simple random walk belongs to $Q_{f(k)}\backslash \partial Q_{f(k)}$ then we set $f(k+1):=f(k)$ and $g(k+1):=g(k)+1$.
\end{quote}

Although this algorithm has an extra randomness due to the SRW, the edges chosen to be revealed in the peeling process are independent of the unknown part, and thus, thanks to Lemma \ref{peelingloi} the process $(|\partial Q_{n}|,|Q_{n}|)_{n\geq 0}$ has the same law as the process obtained with Algorithm  $ \mathcal{L}$. We put  \begin{eqnarray*}\tau_{n} &:=& \sup\{g(k): f(k) = n\}. \end{eqnarray*}
In words, $\tau_n$ is the number of steps made by the SRW  inside $Q_n$. Note that $X_{\tau_n} \in \partial Q_n$.  Since $Q_{n+1}$ differs from $Q_{n}$ by the peeling of an edge of $\partial Q_{n}$ incident to $X_{\tau_{n}}$ we deduce by induction that for every $i \geq 0$ we have\begin{eqnarray} \op{d_{gr}}\left( \partial Q_{i}, \{X_{0}, \ldots , X_{\tau_{i}}\} \right) &\leq& 2.  \label{boundary} \end{eqnarray} 
 

\paragraph{Interpretation.} Let us  recall the definition of the pioneer points of $(X_{n})_{n \geq 0}$. For any $k \geq 0$ we denote by $ \mathcal{R}_{k}$ the submap of $Q_{\infty}$ formed by the faces that are adjacent to $\{X_{0}, X_{1}, \ldots , X_{k}\}$. A moment's thought shows that $ \mathcal{R}_{k}$ is a quadrangulation with holes. Since $Q_{\infty}$ almost surely has one end, only one of these holes corresponds to an infinite quadrangulation and we denote by $ \overline{\mathcal{R}}_{k}$ the quadrangulation obtained from $ \mathcal{R}_{k}$ after filling all the finite holes with their respective quadrangulations in $Q_{\infty}$. Henceforth $ \overline{\mathcal{R}}_{k}$ is a quadrangulation with a simple boundary called the \emph{hull} of the range of $X$ up to time $k$. Recall that for $k \geq 1$, the $k^{th}$ step $X_{k}$ of the SRW is a pioneer point ($k$ is a pioneer time) if   \begin{eqnarray*} X_{k} &\in& \partial \overline{ \mathcal{R}}_{k-1}.  \end{eqnarray*}
 By convention $k=0$ is a pioneer time.
 \begin{proposition} The pioneer times of $(X_{n})_{n\geq 0}$ are exactly the times $\{\tau_{k} : k \geq 0\}.$\end{proposition}
\proof[Proof (Sketch)] We prove by induction the following property : $(*)$ For all $k \geq 0$ such that we have $X_{g(k)}  \in Q_{f(k)} \backslash \partial Q_{f(k)}$ then $Q_{f(k)}$ corresponds to the hull of the range of $X$ up to time $g(k)$. Assume that  this property holds for a certain $k$ and let $n= f(k)$. Obviously, the property holds  for all $k$ such that $g(k) < \tau_{n}$. The time $\tau_{n}$ is thus a pioneer point and the peeling process is then triggered and we discover all the faces in $Q_{\infty} \backslash Q_{n}$ adjacent to $X_{\tau_{n}}$ (and fill the holes they create) until all of $ \overline{ \mathcal{R}}_{\tau_{n}}$ is revealed. At this point the SRW lies inside the current quadrangulation (not on its boundary) and the property $(*)$ holds anew. Details are left to the reader. \endproof 
 
 Notice that the number of peeling steps  is larger than or equal to the number of pioneer points visited so far minus one (recall that $t=0$ is a pioneer time) because the discovery of a pioneer point automatically triggers a new step of peeling. However, the maximal number of steps that a pioneer point can trigger is obviously bounded above by its degree in $Q_{\infty}$, where the degree $ \mathrm{deg}(u)$ of a vertex $u$ is the number of edges adjacent to it. Hence we have  \begin{eqnarray}  n+1 \quad \geq \quad  \#\{ \mbox{pioneer times } \leq \tau_{n} \} & \geq & \frac{n}{  \max \{ \mathrm{deg}(u) : \mathrm{d_{gr}}(\rho, u) \leq n\} }. \label{maxdeg}  \end{eqnarray}
 


\section{Construction of $Q_{\infty}$ from a labeled tree}
\label{schaeffer}
In Section \ref{geo} we gather some geometric estimates on the UIPQ. Most of the results depend on a Schaeffer-like construction of the UIPQ  introduced in \cite{CMMinfini}.  For sake of completeness, we briefly recall it here, the interested reader should consult  \cite{CMMinfini} for more details.

\subsection{The uniform infinite labeled tree $(T_{\infty},\ell)$}We use the standard formalism for plane trees as found in \cite{Nev86}. A plane tree $t$ is a tree given with an ancestor and an order for the children of any vertex $u \in t$. We use the same notation as \cite{CMMinfini}. In particular, the ancestor (or root) of a plane tree $t$ is denoted by $\varnothing$ and its size $|t|$ is its number of vertices. In the following, all the trees that we consider are plane trees.  We denote the set of all rooted plane infinite trees with only one infinite geodesic (also called spine) by $\mathscr{S}$. A tree $t \in \mathscr{S}$ is thus composed of a unique infinite geodesic $$\big\{\varnothing = s_{0},s_{1},s_{2}, \ldots\big\},$$ and finite trees grafted to the left and to the right of each vertex $ s_{i}$.  The degree  of a vertex $u \in t$, denoted by $\op{deg}(u,t)$, is the number of edges adjacent to $u$ in $t$. Such a tree can properly be drawn in plane without accumulation point of the vertices (and matching the ordering of the tree with the clockwise orientation of the plane). A \emph{corner} of a vertex $u \in t$ is an angular sector between two consecutive edges in clockwise order around $u$ (in a plane representation of $t$).  A vertex of degree $k$ thus has  $k$ corners. The contour of the tree $t \in \mathscr{S}$ is the bi-infinite sequence of corners $$\big\{ \ldots , c_{-2},c_{-1},c_{0},c_{1},c_{2},\ldots \big\}$$ sorted in clockwise order where   $c_{0}$ is the corner of the ancestor $\varnothing$ where the tree $t$ is rooted. If $c_{i}$ and $c_{j}$ are two distinct corners of $t$, we denote by $[c_{i},c_{j}]$ the set of corners that are inbetween $c_{i}$ and $c_{j}$ for the contour order (notice that if $i \geq 0$ and $j \leq 0$ then $[c_{i},c_{j}]$ is infinite). If $c$ is a corner of $t$ then $ \mathcal{V}(c)$ is the vertex associated with $c$. 

A \emph{labeling} of a plane tree $t$ is a collection $\{\ell(u) : u \in t\}$ of variables with values in $ \mathbb{Z}$ attached to each vertex of $t$ such that $\ell( \varnothing) =0$ and $\ell(u)-\ell(v) \in \{-1,0,+1\}$ for any neighboring vertices $u,v \in t$. The label $\ell(c)$ of a corner $c$  is the label of its vertex.  We now define the notion of successor. Let $(t,\ell)$ be an infinite labeled tree in $ \mathscr{S}$ and let $c_i$ be a corner of $t$. The successor  of $c_{i}$ is the first corner $ \mathcal{S}(c_{i})$ belonging to $$ \big\{ c_{i+1}, c_{i+2}, \ldots \big\} \cup \big\{ \ldots ,c_{i-2}, c_{i-1} \big\}$$ such  that  $\ell( \mathcal{S}(c_i)) = \ell( c_i)-1.  $ Note that if the vertex associated with $c$ has a minimal label among $t$ then $c$ has no successor (this case will not show up in our setup).

 \bigskip 
We now present the random infinite labeled tree $(T_\infty,\ell)$  which the UIPQ is constructed from. Firstly,  $T_{\infty}$ is a geometric critical Galton-Watson tree conditioned to survive (see \cite{Kes86,LPP95b}). The distribution of $T_{\infty}$ can roughly be described as follows: $T_{\infty}$ has a unique spine (thus $T_{\infty} \in \mathscr{S}$) and the subtrees grafted to the left and to the right of each vertex of the spine are independent critical geometric Galton-Watson trees. See \cite{CMMinfini} for more details.  We recall that if $T$ is a critical geometric Galton-Watson (of parameter $1/2$) then we have 
 \begin{eqnarray} P(|T| = n+1) &=& \frac{ \mathrm{Cat}(n)}{2\cdot 4^n} \quad \sim \quad  \frac{n^{-3/2}}{2 \sqrt{\pi}},   \label{bassin}\end{eqnarray}where $|T|$ is the number of vertices of $T$ and $ \mathrm{Cat}(n) = {2n \choose n}/(n+1)$ is the $n$th Catalan number. 
 Thus the random variable $|T|$ is in the domain of attraction of a completely asymmetric stable variable with parameter $1/2$. 
 
 We then label the tree $T_\infty$ according to the following device: Conditionally on $T_{\infty}$, let $ \{\mathsf{d}_{e} : e \in \op{Edges}(T_{\infty})\}$ be independent random variables uniformly distributed over $\{-1,0,+1\}$ carried by the edges of $T_{\infty}$. This defines a labeled tree $(T_{\infty}, \ell)$ where the label of a vertex  is the sum of the $ \mathsf{d}_{e}$'s  along its ancestral path towards the ancestor $\varnothing$. If $( t,\ell)$ is a labeled tree we set  \begin{eqnarray*} \Delta(t,\ell) &=& \max\{ |\ell(u)| : u \in t\}.  \end{eqnarray*}

  In the following, $(T_{\infty},\ell)$ always denotes the tree constructed above that we call the uniform infinite labeled tree. For $n \geq 0$, we denote by $ \mathcal{T}_{n}$ the (labeled) subtree obtained from $(T_{\infty},\ell)$ after pruning at the $n^{th}$ vertex of the spine $s_{n}$, that is, we remove all the offspring of $s_{n}$ (but we keep $s_{n}$). Recall that $| \mathcal{T}_{n}|$ is the number of vertices of $ \mathcal{T}_{n}$. We also denote  $ \varnothing( \mathcal{T}_{n})  := \max \{  \mathrm{dist}(u,v) : u,v \in \mathcal{T}_{n}\}$ where $ \mathrm{dist}(.,.)$ is the graph distance in $ \mathcal{T}_{n}$, its diameter. Recall the notation $\preceq, \succeq$ and $\approx$ from Section \ref{notations}.
 
 \begin{proposition} \label{arbre}We have
   \begin{eqnarray}  
   \varnothing( \mathcal{T}_{n}) & \approx & n,  \nonumber \\
   | \mathcal{T}_{n}| & \approx & n^2, \label{volume} \\
    \Delta( \mathcal{T}_{n}) & \approx & n^{1/2}. \label{label}
    \end{eqnarray} 
 \end{proposition}
 \proof[Proof (Sketch)] These are pretty standard facts but we include a proof for sake of completeness. 
Let $n \geq 0$. The tree $ \mathcal{T}_{n}$ is composed of the first $n+1$ vertices on the spine together with $2n$ independent critical geometric Galton-Watson trees grafted to the right-hand side and to the left-hand side of $s_{0},s_{1}, \ldots , s_{n-1}$ (when there is no tree on one side of a vertex of spine we consider that we grafted the tree with a single vertex). Thus we have   \begin{eqnarray*} | \mathcal{T}_{n}| = 1-n +\sum_{i=1}^{2n} X_{i} \quad \mbox{ and } \quad n \leq \varnothing( \mathcal{T}_{n})  \leq   2\big(n + \max_{1 \leq i \leq 2n} H_{i}\big)\end{eqnarray*}where $X_{1}, H_{1},X_{2},H_{2,}\ldots$ are respectively the size and the height of the $2n$ critical geometric Galton-Watson trees grafted on the $n$ first vertices of the spine. Recall from \eqref{bassin} that we have 
 $P(X_{1} \geq n) = \Theta(n^{1/2})$ . Recall also the classical Kolmogorov's estimate $P( H_{1} \geq n) \sim n^{-1}$. From the latter we easily deduce using Borel-Cantelli lemma that eventually $H_{i} \leq i \log^{2}(i)$ and thus  $ \varnothing( \mathcal{T}_{n}) \approx n$. Concerning the size $| \mathcal{T}_{n}|$, the analogue of the law of the iterated logarithm in the case of infinite variance (see \cite[Section 3.9]{KK08}) directly show that $S_{n} = \sum_{i=1}^{2n}X_{i} \approx n^2$ which implies $| \mathcal{T}_{n}| \approx n^2$.

 Let us now turn to \eqref{label}. Recall that conditionally on the tree structure of $ \mathcal{T}_{n}$,  the labels evolve along the branches of $ \mathcal{T}_{n}$ as a random walk $(Z_{k})_{k \geq 0}$ whose increments are uniform in $\{-1,0,+1\}$. Looking at the labels of $s_{0}, ... , s_{n-1}$ we deduce that $ \Delta( \mathcal{T}_{n}) \geq \max_{0 \leq i \leq n-1} |\ell( s_{i})|$ which gives the lower bound $ \Delta( \mathcal{T}_{n}) \succeq n^{1/2}$. For the upper bound, we have 
  \begin{eqnarray*} P\Big( \Delta( \mathcal{T}_{n}) > \log^3(n) n^{1/2} \,\left| \,  \mbox{Structure of }{ \mathcal{T}}_{n}\Big) \right. & \leq & |  {\mathcal{T}}_{n}| P \Big(\sup_{0 \leq k \leq \varnothing( \mathcal{T}_{n})} |Z_{k}| \geq \log^3(n) n^{1/2}\Big).  \end{eqnarray*}
On the event $A_{n}:=\{ |\mathcal{T}_{n}| \leq n^{3} \mbox{ and } \varnothing( \mathcal{T}_{n}) \leq n \log^2(n)\}$ the right-hand side of the last display is $O(n^{-2})$.  But  the previous estimates imply that $A_{n}$ eventually occur and thus an application of Borel-Cantelli proves $\Delta( \mathcal{T}_{n}) \preceq n^{1/2}$.    \endproof

 \subsection{Schaeffer construction}
\begin{figure}[!h]
 \begin{center}
 \includegraphics[width=13cm]{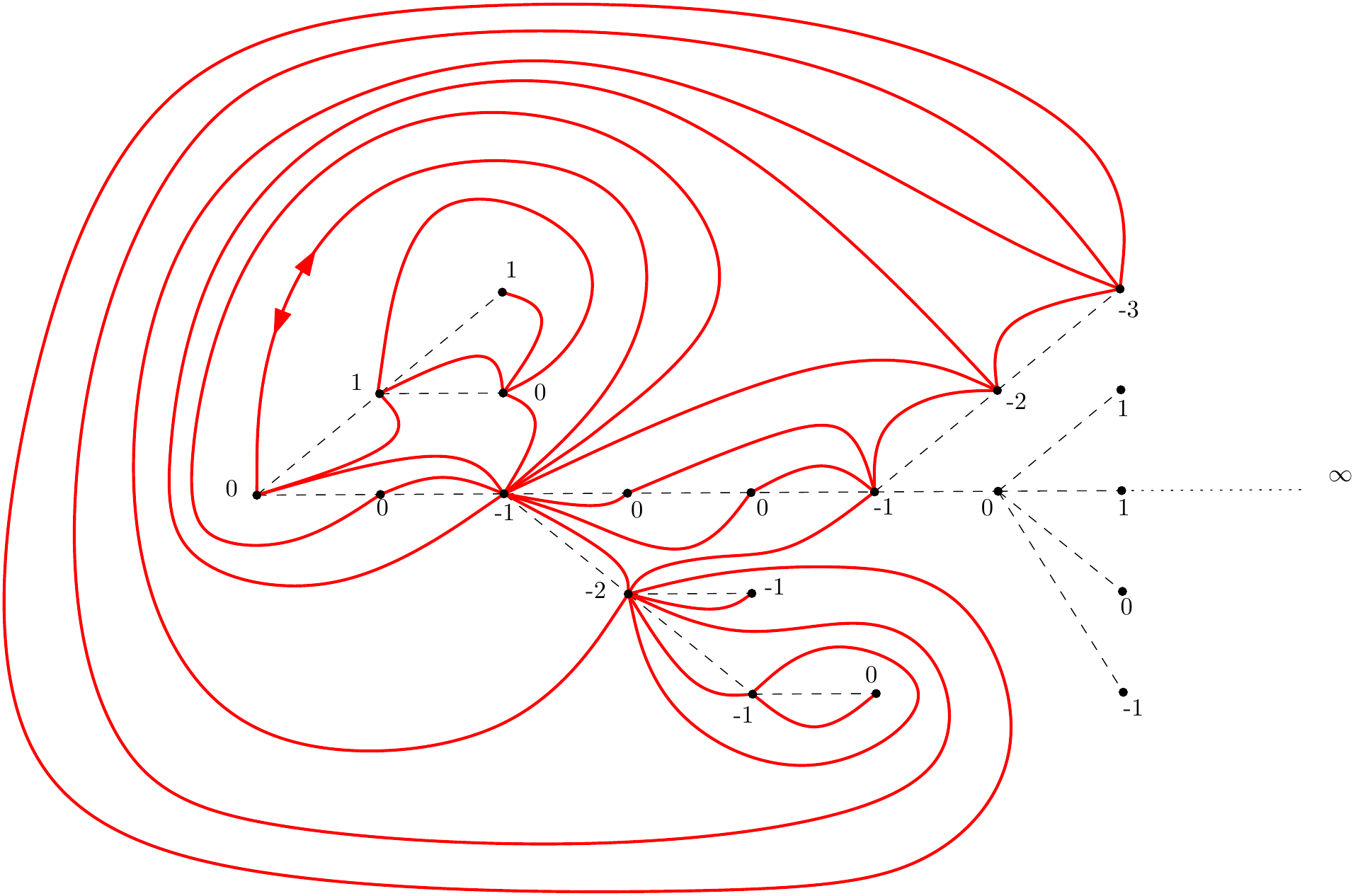}
 \caption{ \label{schaefferfig}Illustration of the construction of  $ \Phi(T_{\infty},\ell)$. The orientation of the root edge is given by an extra Bernoulli variable. }
 \end{center}
 \end{figure}

  A rooted quadrangulation is  associated  with $(T_{\infty}, \ell)$ by the following device. We first embed the labeled tree $ T_{\infty}$ in the plane such that there is no accumulation point for the vertices and such that the edges are not crossing (this is possible since $ {T}_{\infty}$ has one spine).  Then for each corner $c$ of $T_{\infty}$, we draw an edge between $c$ and its successor $ \mathcal{S}(c)$ (note that this successor exists a.s.). All the edges can be drawn in a non-crossing fashion and after erasing the (embedding of the) tree, the resulting map is an infinite quadrangulation. See Fig.\,\ref{schaefferfig}.

 We root it at the edge emanating from the root corner of $T_{\infty}$ whose orientation is given by an extra independent Bernoulli variable $\eta \in \{+,-\}$. The quadrangulation that we obtain, denoted by $\Phi( T_{\infty},\ell)$ (the dependance in $\eta$ is implicit), has the same distribution as $Q_{\infty}$, see \cite{CMMinfini}. In this representation, the vertices of the map are exactly the vertices of the tree $T_{\infty}$, and we shall always make this identification. Using the fact that any neighboring vertices in $ \Phi( T_{	\infty}, \ell)$ must have labels that differ by $1$ in absolute value, we easily get that for every $u,v \in T_{\infty}$ we have  \begin{eqnarray} \label{trivial} \mathrm{d_{gr}}(u,v)  &\geq&  |  \ell(u)-  \ell(v)|.  \end{eqnarray} In fact, the labeling $\ell$ of the vertices of $\Phi(T_{\infty},\ell)$ inherited from this construction has a metric meaning within the quadrangulation $ \Phi(T_{\infty},  \ell)$: The main result of \cite{CMMinfini} states that for every $u,v \in \Phi( T_{ \infty},  \ell)$ we have 
 \begin{eqnarray*}  \ell(u)-  \ell(v) &=& \lim_{ z \to \infty} \big(\mathrm{d_{gr}}(z,u) -  \mathrm{d_{gr}}(z,v)\big),  \end{eqnarray*} where $z \to \infty$ means that $ \mathrm{d_{gr}}_{}(\rho,z) \to \infty$.  We will not use this precise result in what follows, however we will make a great use of the following  bounds on the distances in $\Phi(T_{\infty},\ell)$.  First of all, the very standard bound 
  \begin{eqnarray} \label{sup} \mathrm{d_{gr}}(u,v) & \leq & 2+ \ell(u)+ \ell(v)-2 \max \big\{ \min_{[c,c']} \ell  : \{ \mathcal{V}(c), \mathcal{V}(c')\} = \{u,v\}\big\}, \end{eqnarray} which can be proved as follows. Consider a corner $c_{i}$ of $u$ and a corner $c_{j}$ of $v$ and suppose that $i \leq j$. We construct the path starting from $c_{i}$ and $c_{j}$  following iteratively their successors. These two paths merge at the first corner after $c_{j}$ with label  $\min_{[c_{i},c_{j}]} \ell -1$ and the concatenation of these two paths up to the merging point gives the bound $ \mathrm{d_{gr}}(u,v) \leq 2+ \ell(u)+\ell(v)- 2\min\{ \ell(c) : c \in [c_i,c_j]\}$. The other cases are similar. We also have a lower bound also called \emph{cactus bound}
    \begin{eqnarray}\label{inf} \mathrm{d_{gr}}(u,v) & \geq &   \ell(u) + \ell(v) - 2 \min\big\{\ell(w) : w \in  \llbracket u,v \rrbracket \big\},  \end{eqnarray} where $ \llbracket u,v \rrbracket$ is the geodesic line between $u$ and $v$ in the tree $T_{\infty}$. Let us sketch the idea of the proof of this lower bound, see \cite[Equation (4)]{CMMinfini}. Excluding trivial cases, we consider a vertex $w \in \llbracket u,v \rrbracket \backslash \{u,v\}$ such that $ \ell(w) < \ell(u)$ and $\ell(w) < \ell(v)$. Then choose two corners $c,c'$ of $w$ on both sides of $ \llbracket u,v \rrbracket$. Here also we consider the two paths formed by the successors of $c$ and $c'$. These two paths merge and their concatenation forms a loop separating $u$ from $v$ in $ \Phi(T_{\infty},\ell)$. Thus by Jordan's lemma any path going from $u$ to $v$ in $Q_{\infty}$ must encounter this loop. To finish notice that all the labels on the loop are less than or equal to the label of $w$ and use the bound \eqref{trivial} to conclude. We safely leave the details to the reader.
    \section{Geometric estimates} \label{geo}
    
    In the following, we will consider that the UIPQ is constructed from a uniform infinite labeled tree, that is, we set $Q_{\infty} = \Phi(T_{\infty}, \ell)$. We shall always identify the vertices of $Q_{\infty}
    $ with those of $T_{\infty}$. 
    Unless mentioned, $ \mathrm{d_{gr}}(.,.)$ stands for the graph distance in $Q_{\infty}$. 
\subsection{Uniform estimates on the degrees}
Our first geometric matter concerns the degrees of the vertices in $Q_{\infty}$. Angel \& Schramm proved that the degree of the origin of the UIPT has an exponential tail, see \cite[Lemma 4.1]{AS03}. We shall provide the exact distribution of the degree of the origin of the UIPQ and give a uniform control among all vertices within a given distance from the origin $ \rho$ of $Q_{\infty}$. 

\begin{figure}[!h]
 \begin{center}
 \includegraphics[width=8cm]{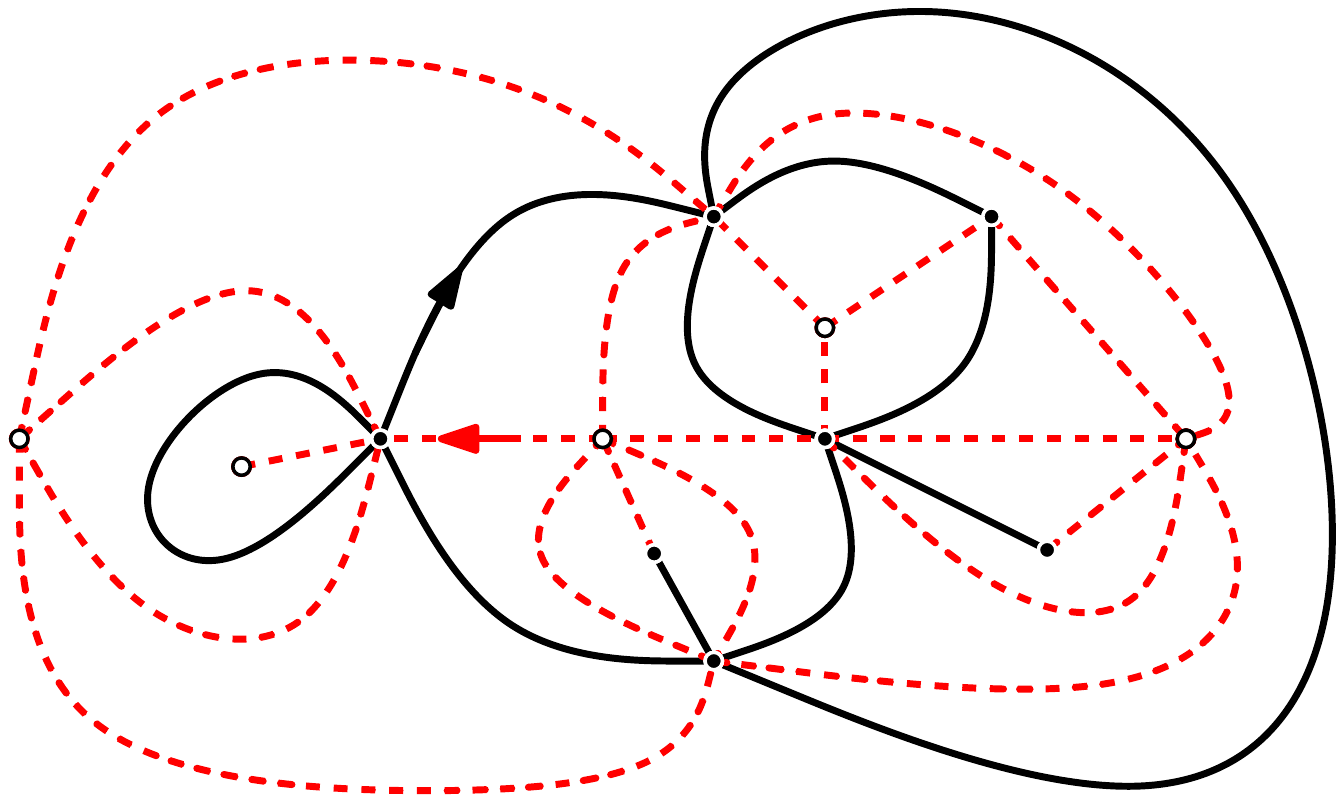}
 \caption{ \label{duality}Illustration of the duality between maps and quadrangulations.}
 \end{center}
 \end{figure}

\begin{proposition} \label{degree} 
\begin{enumerate}[(i)]
\item For every $ y \in (0,6/5)$, we have 
 \begin{eqnarray*} E\Big[y^{\mathrm{deg}(\rho,Q_{\infty})}\Big] &=& 
 \frac{y}{12}(1+y/2)^{-1/2}(1-5y/6)^{-3/2}.  \end{eqnarray*}In particular, $P(\mathrm{deg}(\rho,Q_{\infty})=k) \sim \sqrt{k/40\pi}(5/6)^k$ as $k \to \infty$. \item Furthermore if   $D_{r}$ denotes the maximal degree of a vertex in $ \mathrm{Ball}(Q_{\infty},r)$, then there exists a constant $K_{1}>0$ such that, almost surely
 \begin{eqnarray*}\limsup_{r \to \infty} \frac{D_{r}}{\log(r)} &\leq& K_{1}. \end{eqnarray*}
 \end{enumerate}
\end{proposition}
\proof[Proof of Proposition \ref{degree} part $(i)$] This result  follows from the  enumeration of general planar maps. Indeed, there is a well-known bijection $ \mathcal{D}$ between the set of all rooted planar maps with $n$ edges and the set of all quadrangulations with $n$ faces.  The application $ \mathcal{D}$ can be described as follows: If $m$ is a planar map with $n$ edges, then in each face of $m$ we put an extra point that we link to all the vertices adjacent to this face. We then erase all the edges of $m$ and are left with a quadrangulation $q$ with $n$ faces, see Fig.\,\ref{duality}. 
The root edge of $m$ is the first edge on the right of the rooted edge of $q$ as depicted on Fig.\,\ref{duality}. In this correspondence, the degree of the root face (on the right of the root edge) of $m$ is equal to the degree of the origin of the root edge of $q$. Part $(i)$ then directly follows from \eqref{def:UIPQ} and Theorem 1 of \cite{GR94}. \endproof

 Before going into the proof of Proposition \ref{degree} part $(ii)$ let us give a lemma on $T_{\infty}$.   Recall that the contour of $T_{\infty}$ is denoted by $(\ldots, c_{-2},c_{-1},$ $c_{0},c_{1},c_{2}, \ldots )$ and its spine by $(s_{0}, s_{1}, s_{2}, \ldots )$. Assume that $T_{\infty}$ has been drawn in the plane and consider the sequence of oriented edges $( \ldots,\vec{e}_{-2}, \vec{e}_{-1}, \vec{e}_{0}, \vec{e}_{1}, \vec{e}_{2},\ldots)$ obtained when doing the contour of the tree in clockwise order such that $ \vec{e}_{0}$ is the first oriented edge encountered after the root corner of the tree $T_{\infty}$. Since $T_{\infty}$ almost surely has one spine, for any  oriented edge $ \vec{ \epsilon}$ of $T_{\infty}$ one can say if $ \vec{ \epsilon}$ is poining towards or from infinity, formally if $ \vec{ \epsilon}_{-}$ and $ \vec{\epsilon}_{+}$ are the origin and endpoints of $ \vec{\epsilon}$ the following quantity is well-defined  \begin{eqnarray*} \zeta ( \vec{\epsilon}) &:=& \lim_{z \to \infty} \big( \mathrm{dist}(z, \vec{\epsilon}_{+})- \mathrm{dist}(z, \vec{\epsilon}_{-})\big) \quad \in \quad \{-1,+1\},  \end{eqnarray*} where $ \mathrm{dist}(.,.)$ is the graph distance in $T_{\infty}$ and $z \to \infty$ means that $  \mathrm{dist}(z,\varnothing) \to \infty$. 

\begin{lemma} \label{contour}If $T_{\infty}$ is a critical geometric Galton-Watson tree conditioned to survive then  the variables $\{ \zeta( \vec{e}_{i})\,,\ i \in \mathbb{Z}\}$ are i.i.d.\,Bernoulli variables of parameter $1/2$.\end{lemma}
This lemma  easily follows from \cite[Lemma 4]{CMMinfini} or \cite{LG86}. We leave to the reader the fact that this lemma together with $T_{\infty} \in \mathscr{S}$ completely characterizes the distribution of $T_{\infty}$. In particular, we deduce that for any $k \in \mathbb{Z}$, the tree $T_{\infty}^{(k)}$ consisting of $T_{\infty}$ re-rooted at the corner $c_{k}$ (and the same planar ordering) has the same distribution as $T_{\infty}$,  \begin{eqnarray} \label{reroot}T_{\infty}^{(k)}  & \overset{(d)}{=}& T_{\infty}.  \end{eqnarray}

\proof[Proof of Proposition \ref{degree} $(ii)$]   By part $(i)$, the degree of $\rho$ in $Q_{\infty}$ has an exponential tail. Since $\rho=\varnothing$ with probability $1/2$ we deduce that the degree of $\varnothing$ in $Q_{\infty}$ has an exponential tail as well. By invariance of $T_{\infty}$ under re-rooting, we deduce that there exists some constant $c>0$ such that
 $P\big(\mathrm{deg}( \mathcal{V}(c_{r})) \geq   c\log(r)\big)  =  O(r^{-2}),  $ where $ \mathrm{deg}(.)$ denotes the degree in $Q_{\infty}$. Applying Borel-Cantelli's lemma we deduce that a.s. we eventually have  \begin{eqnarray} \label{test} \mathrm{deg}( \mathcal{V}(c_{i})) &\leq &c \log(r) \quad \mbox{ for all } |i| \leq r.  \end{eqnarray} 
 For $r \geq 1$, let $ \sigma_{r}$ be the first $i \geq 0$ such that the $i^{th}$ vertex along the spine of $T_{\infty}$ has label $\ell(s_{i})=-r$. Recall that the tree $T_{\infty}$ pruned at $s_{ \sigma_{r}}$ is denoted by $\mathcal{T}_{\sigma_r}$. Thanks to \eqref{inf} we deduce that if $v\in Q_{\infty}$ is such that $\op{d}_{\op{gr}}^{}(  \varnothing,v) \leq r-1$ then $v \in \mathcal{T}_{\sigma_r}$. Since the graph distance in $Q_{\infty}$ between $ \varnothing$ and the origin $\rho$ of $Q_{\infty}$ is either $0$ or $1$ we deduce that    \begin{eqnarray} \mathrm{Ball}(Q_{\infty},r)  &\subset&  \{ u \in Q_{\infty} : \mathrm{d_{gr}}(\rho,u) \leq r+1\}  \nonumber \\ &\subset&  \{ u \in Q_{\infty} : \mathrm{d_{gr}}(\varnothing,u) \leq r+2\} \nonumber \\  &\subset &  \mathcal{T}_{\sigma_{r+3}}, \label{bornesup}  \end{eqnarray}  in terms of vertex sets. If $I_{r}$ and $S_{r}$ are the minimal resp.\,maximal indices of a corner belonging to $ \mathcal{T}_{\sigma_r}$ then arguments similar to that of the proof of Proposition \ref{arbre} show that  $S_{r}-I_{r} \approx r^4$ (in fact $S_{r}-I_{r} \preceq r^\alpha$ for some $\alpha >0$ would suffice here). Using this and \eqref{test},  we deduce that there exists a constant $K_{1}$ such that a.s.\,for every $ I_{r} \leq  i \leq S_{r}$ we have $ \mathrm{deg}( \mathcal{V}(c_{i}), Q_{\infty}) \leq K_{1} \log(r)$. Using  \eqref{bornesup}, we complete the proof of the proposition.\endproof

In particular, we deduce from \eqref{maxdeg} and the last proposition that 
 \begin{eqnarray} \label{pioneern} \# \{ \mbox{pioneer times } \leq \tau_{n}\} & \approx & n.  \end{eqnarray}

\subsection{Growth}
 We establish the analogs of the theorem of Angel \cite{Ang03} about the volume growth of the UIPT in the case of the UIPQ.  Recall that $ \mathrm{Ball}(Q_{\infty},r)$ is composed of the faces of $Q_{\infty}$ that have at least one vertex  at distance strictly less than $r$ from the origin $\rho \in Q_{\infty}$, and that $| \mathrm{Ball}(Q_{\infty},r)|
 $ is the number of vertices of $ \mathrm{Ball}(Q_{\infty},r)$. \begin{proposition} \label{growth} We have $| \mathrm{Ball}(Q_\infty,r)| \approx r^4$.  
 \end{proposition}
 See also \cite{CD06,LGM10} for closely related results. 
 \proof Here also we consider that $Q_{\infty}= \Phi(T_{\infty},\ell)$. We begin with the upper bound $ | \mathrm{Ball}(Q_{\infty},r)| \preceq r^4$. Let $r \geq 1$. As in the proof of Proposition \ref{degree} we use the tree $ \mathcal{T}_{\sigma_r}$ consisting of $T_{\infty}$ pruned at the first vertex $s_{ \sigma_{r}}$ of the spine reaching label $-r$. Recall \eqref{bornesup}. 
 Since $ \sigma_{r}$ is the hitting time of $-r$ by a random walk with steps distribution uniform in $\{-1,0,+1\}$, we have $\sigma_{r} = \sigma_{1}^{(1)} + \ldots + \sigma_{1}^{(r)}$ where de $\sigma_{1}^{(i)}$ are i.i.d.\,and distributed as $\sigma_{1}$. Standard calculations show that $P(\sigma_{1} \geq n) \sim C n^{-1/2}$ for some $C>0$. Hence similar arguments as those presented in the proof of Proposition \ref{arbre} show  that $\sigma_{r} \approx r^2$.  We can thus combine this fact together with \eqref{volume} and \eqref{bornesup} to complete the upper bound.

  We now turn to the lower bound. 
  For $ r \geq 1$, 
  we put  \begin{eqnarray*} L_{r}&=& \sup \{i \geq 0 : \Delta(\mathcal{T}_{i}) < r\}. \end{eqnarray*} Consistently we the preceding notation we write $ \mathcal{T}_{L_{r}}$ for the tree $T_{\infty}$ pruned at $s_{L_{r}}$. Using the bound \eqref{sup}, one sees that all the vertices in $ \mathcal{T}_{L_{r}}$ are at a graph distance at most $3r+2$ from $ \varnothing$ in $ Q_{\infty}$, which implies  \begin{eqnarray} \label{borne}  \mathcal{T}_{L_{r}} & \subset &  \mathrm{Ball}(Q_{\infty},3r+4), \end{eqnarray} in terms of vertex sets.
  Using \eqref{label} we deduce that $L_{r} \approx r^2$. Henceforth by \eqref{volume} we have $ | \mathcal{T}_{L_{r}}| \approx r^4$ which together with \eqref{borne} completes the proof of the proposition. \endproof

 \subsection{Tentacles}
Our third estimate deals with the distances in the hull of the ball of radius $r$ in $Q_{\infty}$.  Recall that $ \overline{ \mathrm{Ball}}(Q_{\infty},r)$ is obtained from $ \mathrm{ Ball}(Q_{\infty},r)$ after filling-in all the finite holes. We show that in fact this procedure does not increase the diameter by much, that is, $\overline{\op{Ball}}(Q_{\infty},r)$ does not grow long ``tentacle''.
\begin{proposition} \label{tentacles}  We have $\max \left\{\op{d}_{\op{gr}}( \rho,u) : u \in \overline{\op{Ball}}(Q_{\infty},r) \right\}   \approx r.$
\end{proposition}
\proof  We use the same notation as in the proof of Proposition \ref{growth}. Note that the lower bound $\max \left\{\op{d}_{\op{gr}}( \rho,u) : u \in \overline{\op{Ball}}(Q_{\infty},r) \right\}   \succeq r$ is trivial. For the upper bound, we will strengthen \eqref{bornesup} and prove that    \begin{eqnarray} \label{bornesupplus}\overline{ \mathrm{Ball}}(Q_{\infty},r)  &\subset  &\mathcal{T}_{\sigma_{r+3}},  \end{eqnarray}  in terms of vertex set in $Q_{\infty}$. Indeed consider $s_{\sigma_{r+3}}$ the first vertex on the spine of $T_{\infty}$ with label $-r-3$ and pick $c$ and $c'$ two corners associated with $s_{\sigma_{r+3}}$ from both sides of the spine. We then draw the two paths in $Q_{\infty}$ starting from $c$ and $c'$ by following the chain of successors. These two paths eventually merge. We consider the  cycle $ \mathscr{C}$ formed by the two paths up to the merging point. It is composed of vertices of labels less that $-r-3$ and thus by \eqref{trivial} at distance at least $r+3$ from $\varnothing$. Since $ \mathscr{C}$ separates $Q_{\infty}$ into two parts and  because the minimal graph distance from a point on the cycle to the origin $\rho$ is at least $r+2$, we deduce that $ \overline{ \mathrm{Ball}}(Q_{\infty},r)$ is contained in the finite part  of $Q_{\infty} \backslash \mathscr{C}$ which is included in $ \mathcal{T}_{\sigma_{r+3}}$ in terms of vertex sets.  It follows from \eqref{sup} that   \begin{eqnarray*} \max \left\{\op{d}_{\op{gr}}( \rho,u) : u \in \overline{\op{Ball}}(Q_{\infty},r) \right\} &\leq& 2+3 \Delta( \mathcal{T}_{\sigma_{r+3}}) .\end{eqnarray*} The proof is completed by using \eqref{label} and the fact that $ \sigma_{r}  \approx r^2$. \endproof
  \endproof
  
\begin{rek}  As a  corollary of Proposition \ref{growth} and \ref{tentacles} we have $| \overline{ \mathrm{Ball}}(Q_{\infty},r)| \approx r^{4}$.\end{rek}
  \subsection{Separating cycles}
  Recall the notation $\gamma_{r}$ for the separating cycle ``at distance $r\geq 0$'' from the origin in $Q_{\infty}$ and $|\gamma_{r}|$ for its length. Krikun \cite{Kri05} showed that a slight variant of $|\gamma_{r}|$ is approximately of order $r^2$ and that once renormalized by $r^2$ it converges in distribution towards a $\Gamma(3/2)$ law. Here we use his results to show:
 \begin{proposition}\label{krikun} We have $ |\gamma_{r}| \approx r^2$.
\end{proposition}
\proof In \cite{Kri05}, Krikun studied a separating cycle closely related to our $\gamma_{r}$. More precisely he considered the cycle $\tilde{\gamma}_{r}$ formed by the vertices at distance $r$ from $\rho$ and the diagonals of the faces of $Q_{\infty}$ between them such that $\td_{r}$ separates $\rho$ from the infinite part of the quadrangulation. Since $\gamma_{r}$ and $\tilde{\gamma}_{r}$ are within distance $2$ from each other,  by Proposition \ref{degree}  we have $|\gamma_{r}| \approx | \tilde{\gamma}_{r}|$ and it thus suffices to prove $|\tilde{\gamma}_{r}| \approx r^2$.  Krikun explicitly computed the transition probabilities of $|\tilde{\gamma}_{r}|$ and showed that the process  $(|\tilde{\gamma}_{r}|)_{r\geq 1}$ is a time-reversed critical branching process with offspring distribution in the domain of attraction of a stable distribution of parameter $3/2$. More precisely we have  (Theorem 2 in \cite{Kri05}) 
 \begin{eqnarray} \label{vraikrikun}P\left( |\tilde{\gamma}_{n+r}| = k  \ \Big| \ |\td_r| = l \right) & =& \frac{[t^k] F(t)}{[t^l]F(t)} P \left( \xi_n = l \mid \xi_0=k\right), \end{eqnarray} where $\xi$ is a critical branching process with an explicit offspring distribution  
 and $F(t) = 3/4( \sqrt{(9-t)/(1-t)}-3)$ is the generating function of its stationary measure. In particular, we have \cite[Proof of Corollary 1]{Kri05}
\begin{eqnarray} \label{equation} P\left(  |\td_{r}|=m \right) &  \leq  &  C m^{1/2}r^{-3} \left( 1-\frac{2}{r^2} \right)^m,\end{eqnarray} for some constant $C>0$ (uniform in $m,r$). 
We immediately deduce that  \begin{eqnarray*} P\left(  |\td_{r}| \geq r^2\log(r) \right)  &\leq& Cr^{-3} \sum_{k \geq r^2\log(r)}   k^{1/2}\left(1-\frac{2}{r^2}\right)^k. \end{eqnarray*}
Regrouping the terms in the right hand side by packets of $r^2$ we get for large $r \geq 0$
 \begin{eqnarray*} {P}\left(  |\td_{r}| \geq r^2\log(r) \right) &\leq& C' \sum_{k \geq \log(r)}k^{1/2} e^{-2k} \quad = \quad   O(r^{-2+\varepsilon}).  \end{eqnarray*}  A direct application of Borel-Cantelli's lemma shows that $ |\td_{r}| \leq \log(r) r^2$ eventually which proves the upper bound of the proposition. The lower bound is a bit more involved. First of all, it follows from \eqref{equation} that ${P}( |\td_{r}| \leq r^2\log^{-1}(r)) = O(\log^{-3/2}(r))$. Applying Borel-Cantelli's lemma  we get that eventually $|\td_{r}| \geq r^2\log^{-1}(r)$ along the values of $r = 2,4,8,\ldots , 2^k, \ldots $. Notice however that the random process $(|\td_{r}|)_{r\geq 1}$ is \emph{not} increasing and thus we cannot interpolate between values of $2^k$. We bypass this problem by using the Markovian nature of $(\tilde{\gamma})$.
 
 To simplify notation we set $l(r) = \lfloor r^2 \log^{-1}(r) \rfloor $ and $u(r) = \lceil r^2 \log(r) \rceil$. We just proved that a.s. we eventually have   \begin{eqnarray}  l(2^k) \quad \leq \quad  |\td_{2^k}| \quad \leq \quad u(2^k).   \label{interpoltat}\end{eqnarray} For $r \geq 1$ we let $A_{r}$ be the following event   \begin{eqnarray*} A_{r} := \left\{ \begin{array}{l} l(r) \leq |\td_{r}| \leq u(r),\\ 
 |\td_{i}| \leq i^2 \log^{-10}(i), \ \mbox{ for some } r\leq i \leq 2r,\\ l(2r) \leq |\td_{2r}| \leq u(2r). \end{array}\right\}.  \end{eqnarray*} 
We claim that $P(A_{r}) = O( \log^{-2}(r))$. This is sufficent to finish the proof of Proposition \ref{krikun}: By applying Borel-Cantelli's lemma to the sequence of events $A_{2^k}$ for $k =1,2,3, \ldots$ we deduce that $A_{r}^c$ eventually holds which combined with \eqref{interpoltat} yields to $|\tilde{\gamma}_{r}| \approx r^2$. Let us now prove the claim. By definition, $P(A_{r})$ is equal to 
 \begin{eqnarray*} &=& \sum_{a = l(r) }^{u(r)}\sum_{b = l(2r) }^{u(2r)}P( |\td_{r}| =a) P\left(\exists r \leq i\leq 2r : |\td_{i}| \leq r^2 \log^{-10}(r)  \ \mbox{and} \  |\td_{2r}|=b \ \Big | \  |\td_{r}| =a \right),\\ 
 &=&  \sum_{a = l(r)  }^{u(r) }\sum_{b = l(2r) }^{ u(2r)}\frac{[t^b] F(t)}{[t^a]F(t)} P( |\td_{r} |=a) P\left(\exists 0 \leq i \leq r :\xi_{i} \leq r^2 \log^{-10}(r)  \ \mbox{and} \ \xi_{r}=a \ \Big | \ \xi_{0} =b \right),  \end{eqnarray*} by \eqref{vraikrikun}. Using standard singularity analysis one shows that $[t^m]F(t) \sim 3/\sqrt{2\pi} m^{-1/2}$ as $m \to \infty$. Using this and \eqref{equation} we can  bound the term $[t^b]F(t)/[t^a]F(t) P( |\td_{r}|=a)$ in the last display by $ C r^{-2}\log^{3/2}(r)$ for some constant $C>0$ uniform in $r \geq 1$. Thus we have 
  \begin{eqnarray*} P(A_{r}) & \leq & C r^{-2} \log^{3/2}(r) \sum_{b = \lfloor l(2r) \rfloor }^{\lfloor u(2r) \rfloor} P\big(\exists 0 \leq i\leq r : \xi_{i} \leq r^2 \log^{-10}(r), l(r) \leq \xi_{r} \leq u(r) \mid \xi_{0} =b \big).  \end{eqnarray*}
  Fix $ l(2r) \leq b \leq u(2r)$ and let us estimate the probability that the branching process $\xi$ starting from $\xi_{0}=b$ reaches a level lower than $r^2 \log^{-10}(r)$ for some $0 \leq i \leq r$ and finally ends at a state $ l(r) \leq \xi_{r}$. Since $\xi$ is a critical branching process, it is in particular a martingale. Thus if we introduce the stopping times  $T_{r} = \inf \{n \geq 0 : \xi_{n} \geq l(r)\}$ and $\tau = \inf\{n \geq 0 : \xi_{n} =0 \}$ we deduce that 
   \begin{eqnarray*} P( T_{r} < \tau < r \mid \xi_{0} = i)  & \leq & \frac{i}{l(r)}.\end{eqnarray*}
   As a consequence, applying the Markov property of $\xi$ at the first time $j$ where $\xi_{j} \leq r^2 \log^{-10}(r)$ we deduce that the event $ \{\xi_0=b \to \xi_j \leq r^2 \log^{-10}(r) \to \xi_r \geq l(r)\}$ for the branching process $\xi$ has a probability less than or equal to $r^2 \log^{-10}(r) / l(r)  \sim \log^{-9}(r)$. Gathering-up the pieces we finally get that $P(A_{r}) = O (\log^{-2}(r))$ as desired.  \endproof

\subsection{Aperture after peeling}
Let $Q_{0} \subset Q_{1} \subset Q_{2} \subset \ldots$ be the sequence of quadrangulations with boundary obtained by a peeling of $Q_{\infty}$. For each $n \geq 0$ we know that $Q_{\infty} \backslash Q_{n}$ has the same distribution as a UIPQ of the $ |\partial Q_{n}|$-gon.  Following \cite{CMboundary}, if $q$ is a quadrangulation with a boundary, we denote the maximal distance between any pair of points of $\partial q$ by $\op{aper}(q)$ and call it \emph{the aperture} of $q$. 
\begin{proposition}\label{aperture}  We have $ \mathrm{aper}( Q_{\infty}\backslash Q_{n}) \preceq | \partial Q_{n}|^{1/2}.$
\end{proposition}

We slightly abuse notation in the last proposition. Of course the reader would have understood that $\mathrm{aper}( Q_{\infty}\backslash Q_{n}) \preceq | \partial Q_{n}|^{1/2}$ means that there exists a constant $\kappa>0$ such that almost surely we eventually have $\mathrm{aper}( Q_{\infty}\backslash Q_{n})  \leq \log^\kappa (n) | \partial Q_{n}|^{1/2}$. 

 \proof We already recalled that for every $n\geq 0$ the quadrangulation with a simple boundary $Q_{\infty} \backslash Q_{n}$ has the same distribution as a UIPQ of the $ |\partial Q_{n}|$-gon. We now recall an estimate of  \cite{CMboundary}: 
\begin{theorem}[{\cite{CMboundary}}]
There exists $c,c'>0$ such that for all $p\geq 1$ and $\lambda >0$ the  aperture of a uniform infinite planar  quadrangulation with simple boundary of perimeter $2p$ satisfies
\begin{eqnarray*} P\left( \op{aper}({Q}_{\infty,2p}) \geq  \lambda \sqrt{p} \right)  &\leq& c \,  p^{2/3}\, \exp\big(-c'\lambda^{2/3}\big) . \end{eqnarray*}
\end{theorem}Thus taking $\lambda = \log^3(n)$ in this theorem and noticing that $|\partial Q_{n}|$ is deterministically less that $2n$, an application of Borel-Cantelli's lemma finishes the proof.	\endproof 


 \section{Remaining proofs} \label{proofs}
We begin with the proof of Theorem	\ref{estimpeeling}.
\subsection{Peeling estimate}
\proof[Proof of Theorem \ref{estimpeeling}] Because of Lemma \ref{peelingloi} it is sufficient to prove  Theorem  \ref{estimpeeling} for \emph{one} peeling algorithm. We thus consider $Q_{0}, Q_{1}, \ldots$ the peeling of $Q_{\infty}$ using Algorithm $ \mathcal{L}$ of Section \ref{peelL}. During this peeling we know from Proposition \ref{ouff} that all the edges on the separating cycles $\{\gamma_{r}, r \geq 0 \}$ must be part of the boundary of some $Q_{n}$. Furthermore, for all $n\geq 0$, any edge on the boundary of $Q_{n}$ must be at a graph distance less that $2$ from a separating cycle $\gamma_{r}$ for some $r \geq 0$. Using the uniform estimates on the degree (Proposition \ref{degree}) and Proposition \ref{ouff} (and the remark after it) we deduce that after $n$ steps of peeling, the boundary $\partial Q_{n}$ of $Q_{n}$ is located at a graph distance less than $2$ from some $\gamma_{R_{n}}$ with  
 \begin{eqnarray*} \sum_{i=1}^{R_{n}} |\gamma_{i}| & \approx & n.  \end{eqnarray*} Coupling the last display with the fact that $|\gamma_{i}| \approx i^2$ (Proposition \ref{krikun}) we deduce that $R_{n} \approx n^{1/3}$.\\
 Using this with Proposition \ref{krikun} and \ref{degree} again, we deduce that since $ \mathrm{d_{gr}}( \partial Q_{n}, \gamma_{R_{n}}) \leq 2$ we have $ |\partial Q_{n}| \approx n^{2/3}$. This proves the first half of Theorem \ref{estimpeeling}.\\
Let us now focus on the volume of $Q_{n}$. From the deductions made above we have $$ |\overline{ \mathrm{Ball}}(Q_{\infty},R_{n}-3)| \leq |Q_{n}| \leq |\overline{ \mathrm{Ball}}(Q_{\infty},R_{n}+3)|.$$ We then use $R_{n} \approx n^{1/3}$ and the remark after Proposition \ref{tentacles} to get  $|Q_{n}| \approx n^{4/3}$. This completes the proof of Theorem \ref{estimpeeling}. \endproof
\subsection{Pioneer points and subdiffusivity}	
With all the estimates that we now have in our hands, the proof of Theorem \ref{main} is effortless. 
\proof[Proof of Theorem \ref{main}] Let $Q_{\infty}$ be the uniform infinite planar quadrangulation and conditionally on $Q_{\infty}$, let $(X_{n})_{n\geq 0}$ be a nearest-neighbor simple random walk starting from the origin $ \rho \in Q_{\infty}$. We consider $Q_{0}, Q_{1}, \ldots$ the peeling $Q_\infty$ according to  Algorithm $ \mathcal{W}$ of Section \ref{peelW}. By Theorem \ref{estimpeeling} we have $|\partial Q_{n}| \approx n^{2/3}$ and applying Proposition \ref{aperture} we deduce that $ \mathrm{aper}(Q_{\infty}\backslash Q_{n}) \preceq n^{1/3}$. Let $D^-_{n}$ and $D^+_{n}$ be the minimal and maximal distance to the origin $\rho \in Q_{\infty}$ of a vertex in $\partial Q_{n}$. Since we have $ \mathrm{aper}(Q_{\infty}\backslash Q_{n}) \geq D_{n}^+-D_{n}^-$ we deduce that $D^+_{n}-D_{n}^- \preceq n^{1/3}$. From the inclusions 
  \begin{eqnarray*} \overline{ \mathrm{Ball}}(Q_{\infty},D_{n}^- -1 ) \quad \subset \quad Q_{n} \quad \subset\quad  \overline{ \mathrm{Ball}}(Q_{\infty},D_n^++1),  \end{eqnarray*} and since $|Q_{n}| \approx n^{4/3}$ (by Theorem \ref{estimpeeling}) we deduce using Proposition \ref{growth} and \ref{tentacles} that $ D_{n}^{-} \preceq n^{1/3}$ and $ D_{n}^+ \succeq n^{1/3}$. But since $ D_{n}^+-D_{n}^- \preceq n^{1/3}$ we must have $D_{n}^+ \preceq n^{1/3}$. To finish the proof, just recall that $Q_{n}$ is the quadrangulation discovered when $n$ steps of peeling have been demanded and that during this time the SRW has discovered $\approx n$ pioneer points by \eqref{pioneern}. Hence the pioneer points discovered so far are contained in $ \overline{ \mathrm{Ball}}(Q_{\infty},D_{n}^+ +1 )$ which has a diameter $\approx n^{1/3}$ by Proposition \ref{tentacles}. Finally, at least one pioneer points is at distance at least $D_{n}^+-2$ from $\rho$, hence  \begin{eqnarray*} \max_{1 \leq i \leq n} \mathrm{d_{gr}}(\rho,P_{i}) & \approx & n^{1/3}.  \end{eqnarray*}
   \endproof

\begin{figure}[!h]
 \begin{center}
 \includegraphics[width=10cm]{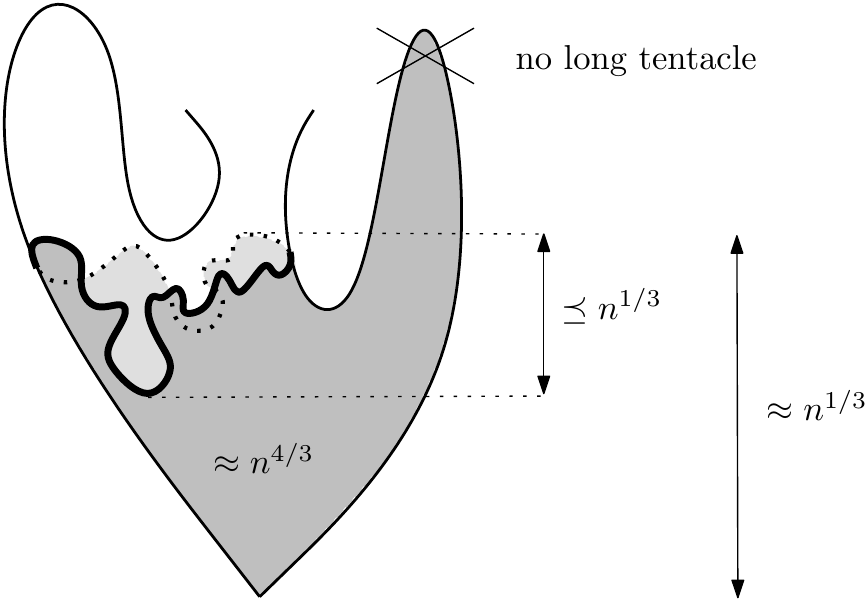}
 \caption{Illustration of the proof of Theorem \ref{main}. The curvy line represents the boundary of $Q_{n}$.}
 \end{center}
 \end{figure}
 
 \proof[Proof of Corollary \ref{subdiffusive}] With the notation of the proof of Theorem \ref{main} the range $\{X_{1}, X_{2}, \ldots , X_{\tau_{n}}\}$ is contained in $ \overline{\mathrm{Ball}}(Q_{\infty}, D_{n}^++1)$. We obviously have 
$ \tau_{n} \succeq n$ (and in fact $\tau_{n}$ could be much larger than $n$). We then use $ D_{n}^+ \approx n^{1/3}$ and  Proposition \ref{tentacles} to conclude. \endproof
 
 \begin{rek} It is clear from the proof of Corollary \ref{subdiffusive} that the $1/3$ exponent of subdiffusivity is not likely to be sharp. Indeed, most of the times are not pioneer times for the simple random walk since two pioneer times could be separated by a long period of time. Yet this phenomenon is hard to control. 
 \end{rek}

\section{Comments and questions}

Before making a more precise list of comments, let us emphase the fact that we focused on the UIPQ for sake of simplicity and because many tools are already available for this model. There should not be major conceptual problem in generalizing our result to other type of random lattices such as the UIPT or Boltzmann maps -- but the required technics might be (much!) more difficult to work with. 
\label{comments}
\subsection{Peeling} The proof of Theorem \ref{main} is not specific to the peeling with algorithm $ \mathcal{W}$ and can be generalized to show that for any peeling  $Q_{0} \subset Q_{1} \subset Q_{2} \subset \ldots$ of $Q_{\infty}$  we actually have
 \begin{eqnarray*} \max_{u \in Q_{n}} \mathrm{d_{gr}}(\rho,u) & \approx & n^{1/3}.  \end{eqnarray*}
 In particular, this result can be applied with other peeling procedures among which:
\begin{itemize}
\item The peeling along layers of $Q_{\infty}$ (giving back a few estimates of Section \ref{geo}),
\item The peeling along a percolation interface as developed in \cite{ACpercopeel},

\item The peeling associated with internal diffusion limited aggregation on the UIPQ,
\item The peeling along a Brownian motion on the Riemann surface associated with $Q_{\infty}$, see \cite{GR10},
\item \ldots
\end{itemize}

\paragraph{Limit processes.} In Theorem \ref{estimpeeling} we established the rough estimates $ |\partial Q_{n}| \approx n^{2/3}$ and $|Q_{n}| \approx n^{4/3}$. One can ask for a precise limit theorem of the re-normalized processes $$ \big(n^{-2/3} |\partial Q_{\lfloor	n t\rfloor}|,n^{-4/3} |Q_{\lfloor	n t\rfloor}|\big)_{ t\geq 0}.$$  Note that the two components are not independent and that Angel \cite{Ang03} conjectured that (in the triangulation case) the first component converges towards a stable process of parameter $3/2$ conditioned to remain positive see \cite{Ber06}. 

\paragraph{Greedy peeling.} Another useful property that has to be addressed about the peeling process is the following. For any peeling $Q_{0} \subset Q_{1} \subset Q_{2}\subset  \ldots$  of $Q_{\infty}$ show that we have 
 \begin{eqnarray*} \bigcup_{n\geq 0} Q_{n} & = & Q_{\infty}.  \end{eqnarray*}
 In words, whatever the algorithm used to peel $Q_{\infty}$, we eventually discover the whole quadrangulation $Q_{\infty}$. This would be implied by the fact that the during the peeling there exist infinitely many times such that $| \partial Q_{n+1}| \leq | \partial Q_{n}|/2$. This result would have  nice applications: Applying it with Algorithm $ \mathcal{W}$ it should  imply that  the range of a simple random walk $(X_{n})_{n\geq 0}$ creates infinitely many loops separating the origin $\rho \in Q_{\infty}$ from $\infty$ a.s.. In particular, two independent simple random walk paths on $Q_{\infty}$ would intersect showing that $Q_{\infty}$ is almost surely Liouville, see \cite{BCstationary}. We expect a similar result to hold for the range of Brownian motions on the Riemann surface of $Q_{\infty}$ thus yielding a different perspective on the result of \cite{GR10}. We hope to pursue these goals in future works.


\subsection{Sudiffusivity}
Note that our subdiffusivity result (Corollary \ref{subdiffusive}) was not based on resistance nor heat kernel estimates as it is generally the case. In reward we can give bounds on the probability that a simple random walk returns to the origin in $n$ steps. For any $x,y$ in $Q_{\infty}$ and $n\geq0$, we denote by $p(x,y,n)$ the probability that a SRW started at $x$ hits the point $y$ at time $n$. Note that $p(x,y,n)$ is random. Our main result implies, 
\begin{cor} We have $ p(\rho,\rho,2n) \succeq n^{-4/3}$.
\end{cor}
\proof 
Recall that $D_{r}$ denotes the maximal degree within distance $r$ of $ \rho$. For any $r \geq 1$ we have 
\begin{eqnarray}
p( \rho, \rho,2n) & \geq & \sum_{x \in \op{Ball}(Q_{\infty},r)} p( \rho, x,n)p(x, \rho,n) \nonumber \\
 &=& \sum_{x \in \op{Ball}(Q_{\infty},r)}  \frac{\op{deg}( \rho)}{\op{deg}(x)}p( \rho,x,n)^2 \nonumber  \\ &\geq  &  {D_{r}^{-1}} \left(\sum_{x \in \op{Ball}(Q_{\infty},r)} p( \rho, x, n) \right)^2  |\op{Ball}(Q_{\infty},r)|^{-1} \nonumber \\
 &\geq & \frac{P(X_{n} \in \mathrm{Ball}(Q_{\infty},r))^2}{D_{r} | \mathrm{Ball}(Q_{\infty},r)|},\end{eqnarray} where we used Cauchy-Schwarz inequality to go from the second to the third line.  Taking $ r =  \lfloor n^{1/3} 
 \log^\kappa(n)\rfloor$ for some $\kappa >0$ we deduce from Proposition \ref{degree}, Proposition \ref{growth} and Corollary \ref{subdiffusive} that $p (\rho, \rho,2n)$ is asymptotically larger than $ n^{-4/3} \log^{\kappa'}(n)$ for some $\kappa' >0$. This completes the proof of the corollary.
\endproof

\begin{rek} Notice that we can produce a lower bound on the displacement of the SRW on the UIPQ by using the crude fact that the electrical resistance $R_{x,y}$ between two points $x,y \in Q_{\infty}$ is less than or equal to $ \mathrm{d_{gr}}(x,y)$, see also \cite{Bar04}. Let us for example given an upper bound on the mean of   \begin{eqnarray*} E_r &=& \inf\{ n \geq 0 : X_n \notin \mathrm{Ball}(Q_\infty,r)\}.  \end{eqnarray*} 
By the result of \cite{CRRST96} we have  \begin{eqnarray*} E[T_{r}]  &\leq& 2 | \mathrm{Ball}(Q_{\infty},r)| R_{\rho, \gamma_{r}} \quad \leq \quad r^{5+o(1)}.  \end{eqnarray*} 
We do not sharpen this result because we do not believe that this is the right  exponent.
\end{rek}

The subdiffusive behavior of the SRW on $Q_{\infty}$ established in Corollary \ref{subdiffusive} is not sufficient to conclude recurrence of the walk. Still, we believe that $Q_{\infty}$ \emph{is} recurrent (see Conjecture \ref{recurrence}) and that the subdiffusivity exponent is critical for deciding recurrence or transience, see Conjecture \ref{1/4}.

We also suspect that one does not need the detailed structure of the UIPQ to establish subdiffusivity but only the existence of bottlenecks at all scales. In particular, is it the case that any planar stationary random graph (see \cite{BCstationary}) with volume growth bigger than quadratic is subdiffusive for the simple random walk? See related conjectures in \cite{BP11}.


\subsection{KPZ}
This part is heuristic. For a mathematically precise statement of the KPZ relations, the reader should consult \cite{DS09}.
\paragraph{Verification of KPZ relation for pioneer exponents.} 

The famous KPZ relation \cite{KPZ88} predicts that certain exponents of statistical mechanics models on a random planar map are related to the analogous exponents on a regular lattice, see \cite{Dup06}. More precisely, let $F$ be a random fractal on a Euclidean space  (for example the set of pioneer points of a Brownian motion). If $F$ has ``dimension'' $2(1-x)$ that means, roughly speaking, that $\varepsilon^{-2(1-x)}$ balls of radius $ \varepsilon$ are necessary to cover $F$ when $ \varepsilon \to 0$. 
Then $x$ is called the \emph{Euclidean scaling exponent of $F$} \cite{DS09}. Similarly, if we consider the same random fractal on a random geometry one can define its ``quantum scaling exponent'' to be $\Delta$ if the number of balls of radius $\varepsilon$ (in the random geometry) needed to cover $F$ is approximatively $(n_{\varepsilon})^{(1-\Delta)}$ where $n_{\varepsilon}$ is the number of balls needed to cover the full space.  The KPZ relation then predicts   \begin{eqnarray*} x &=& \frac{\gamma^2}{4} \Delta^2 + \left(1- \frac{\gamma^2}{4}\right) \Delta,  \qquad \qquad \qquad  \mbox{(KPZ)}\end{eqnarray*} where $0 \leq \gamma<2$ is a parameter depending on the features of the model that produced the random fractal. In particular, in the case of fractals coming from a Brownian motion  we should have $\gamma = \sqrt{8/3}$.  

Going to a discrete level, a random subset $ \mathscr{F}_{n}$ of a planar quadrangulation with $n$ faces is said to have a quantum scaling exponent $ \Delta_{D}$ if $| \mathscr{F}_{n}|$ is of order $ n^{1-\Delta_{D}}$ as $n \to \infty$. Taking a ball of radius $r$ in the UIPQ, we know by Theorem \ref{main} that $ \approx r^3$ pioneer points are visited before the walk exits this ball which contains $\approx r^4$ points. Putting this together, we deduce that the discrete quantum scaling exponent for pioneer points is $\Delta_{D} = 1/4$. Going through (KPZ) this becomes $x_{D} = 1/8$. Indeed, $2-2 x_{D} = 7/4$ is the dimension of the set of pioneer points of the Brownian motion as identified by Lawler, Schramm, Werner  \cite{LSW01} in the Euclidean case. Notice that various quantum scaling exponent for simple random walk on random lattices were derived non-rigorously by Duplantier \& Kwon \cite{DK88}.

\paragraph{Support for Conjecture \ref{1/4}.} Let us use once more the KPZ relation for intersection exponents of simple random walks. More precisely, the probability that $L\geq 1$ independent random walks starting from the same point in a random lattice are not intersecting each other up to time $n$ is supposed to decay as $ n^{- \Delta_L+o(1)}$. These exponents can be derived from the Euclidean case \cite{LSW01} using the KPZ relation and we have \cite[Eq (3.14)]{Dup06}
 \begin{eqnarray*} \Delta_L &=& \frac{1}{2}\left(L-\frac{1}{2}\right).  \end{eqnarray*}  The special case $ \Delta_1 = 1/4$ corresponds to the disconnection exponent, meaning that the probability that the origin of \emph{one} walk has not been disconnected from infinity after $n$ steps decays as $n^{-1/4+o(1)}$. By time reversing this propability is also that of the $n^{th}$  step of the walk being a pioneer point. Thus we should have  \begin{eqnarray*} {P}(n \mbox{ is a pioneer time}) &\asymp& n^{-1/4 +o(1)}.  \end{eqnarray*} Henceforth, in $n^4$ steps the SRW should have discovered roughly $\sum_{k=0}^{n^4} k^{-1/4} \asymp n^3$ pioneer points and by Theorem \ref{main} the maximal displacement from the root in the first $n^4$ steps is $\approx n$. This supports  Conjecture \ref{1/4}.  Note  that this is also equivalent to the fact that the KPZ relation sends $x=0$ to $\Delta=0$, in other words, if the simple random walk covers most of the lattice in the Euclidean case (say for example covers most of the ball of radius $r^{1-o(1)}$ before exiting the ball of radius $r$) then it should be the same in the random lattice case. 

\paragraph{SAW $vs$ SRW.} One of the keys to our result is that after discovering a certain part $Q_{n}$ of the UIPQ which corresponds to the hull of a simple random walk (but could be the hull of a percolation cluster \ldots) then the unknown quadrangulation with a boundary $Q_{\infty}\backslash Q_{n}$ is independent of $Q_{n}$ conditionally on the length of the boundary. 

The boundary of $ \partial (Q_{\infty} \backslash Q_{n})$ can be seen as a self-avoiding loop surrounding the discovered part. Indeed, the annealed model of self-avoiding walk (SAW) on a quadrangulation is totally equivalent to the model of quadrangulation with a simple boundary, just zip the boundary or unzip the SAW, see \cite{CMboundary}. Heuristically speaking, we see that locally the boundary of the range of a simple random walk on the UIPQ is, in a certain sense,  close to a self-avoiding walk. This fact is conjectured in planar Euclidean geometry, but still open.


\endproof

\bibliographystyle{alpha}

\begin {tabular}{l p{3mm} l}
Mathematics Department && D\'epartement de Math\'ematiques et Applications \\
The Weizmann Institute  && Ecole Normale Sup\'erieure, 45 rue d'Ulm \\
Rehovot 76100, Israel &&  75230 Paris cedex 05, France
\end {tabular}

\medbreak 
\noindent nicolas.curien@ens.fr\\
\noindent itai@wisdom.weizmann.ac.il\\

\end{document}